\documentclass[11pt, final]{amsart}
\usepackage{amssymb}
\usepackage[mathscr]{eucal}
\usepackage{amscd}
\usepackage{amsmath}
\usepackage{cite}
\usepackage{setspace}
\usepackage{tikz}
\usetikzlibrary{shapes,arrows,fit,calc,positioning}
\tikzset{box/.style={draw, rectangle, thick, text centered, minimum height=3em}}
\tikzset{line/.style={draw, thick, -latex'}}

\def\Ker{\operatorname{Ker}}

\def\dim{\operatorname{dim}}
\def\supp{\operatorname{supp}}

\theoremstyle{plain}
\newtheorem{theorem}{Theorem}

\newtheorem{proposition}{Proposition}

\newtheorem*{psp}{The Principle of Small Perturbations}
\newtheorem{lemma}{Lemma}
\newtheorem{corollary}{Corollary}
\newtheorem{definition}{Definition}
\theoremstyle{remark}

\newtheorem{example}{Example}

\newcommand{\jinf}{_{j=1}^{\infty}}

\newcommand{\lims}{\lim_{ \substack{ \nu \to \infty \\ \nu<k_1<...<k_n}}}

\title{The Kalton-Peck space as a spreading model}

\address{Departamento de Matem\'aticas, Universidad de Extremadura, Avenida de Elvas s/n, 06011 Badajoz, Spain.}
\email{jesus@unex.es}
\author[J. Su\'arez]{Jes\'us Su\'arez}
\thanks{The author was supported by no institution.}
\begin{document}
\maketitle
\begin{center}
\textit{To Zulema, who is always ready to smile}
\end{center}
\begin{abstract} The so-called Kalton-Peck space $Z_2$ is a twisted Hilbert space induced, using complex interpolation, by $c_0$ or  $\ell_p$ for any $1\leq p\neq 2<\infty$. Kalton and Peck developed a scheme of results for $Z_2$ showing that it is a very rigid space. For example, every normalized basic sequence in $Z_2$ contains a subsequence which is equivalent to either the Hilbert copy $\ell_2$ or the Orlicz space $\ell_M$. Recently, new examples of twisted Hilbert spaces, which are induced by asymptotic $\ell_p$-spaces, have appeared on the stage. Thus, our aim is to extend the Kalton-Peck theory of $Z_2$ to twisted Hilbert spaces $Z(X)$ induced by asymptotic $c_0$ or $\ell_p$-spaces $X$ for $1\leq p<\infty$. One of the novelties is to use spreading models to gain information on the  isomorphic structure of the subspaces of a twisted Hilbert space. As a sample of our results, the only spreading models of $Z(X)$ are $\ell_2$ and $\ell_M$, whenever $X$ is as above and $p\neq 2$.

\end{abstract}

\tableofcontents
\newpage
\section{Introduction}
The topic of twisted Hilbert spaces, i.e. Banach spaces $Z$ containing an isomorphic copy of the Hilbert space $\ell_2$ such that the quotient by this copy $Z/\ell_2$ is isomorphic to $\ell_2$, is well established within geometric nonlinear functional analysis. The book of Benyamini and Lindenstrauss contains the basic main points \cite[Chapter 16]{BenLin}.

The story began with the solution of the Palais problem with the Enflo, Lindenstrauss and Pisier's example \cite{ELP} and later on continued with the irruption of the Kalton-Peck space $Z_2$ in \cite{KaPe}. We would like to recall that the Kalton-Peck space is a major achievement, root of many counterexamples and it must appear in any catalogue of the indispensable Banach spaces for any researcher.

Kalton subsequently provided us in his monumental \cite{ka} with a general procedure to construct twisted Hilbert spaces using complex interpolation. Roughly, if $X$ denotes a separable Banach space for which complex interpolation yields a Hilbert space in the form $(X, X^*)_{1/2} = \ell_2$, then a twisted Hilbert space arises as the derivation of the previous formula. It is known as the derived space of the interpolation couple $(X,X^*)$ at $1/2$ and denoted $Z(X)$. Formally, the derivative gives rise to a (usually nonlinear) map called a centralizer and denoted $\Omega$ which governs the norm of $Z(X)$. Thus, such $Z(X)$ may be defined as the set of couples $(x,y)$ such that $$\|x-\Omega(y)\|_{\ell_2}+\|y\|_{\ell_2}<\infty,$$ 
where $x,y\in \ell_2$. Summing up, as well as Hilbert spaces are important, then so does twisted Hilbert spaces since they appear as a derivation of a Hilbert space.

Some papers have appeared recently relaunching the topic with new examples and in new directions, see for example \cite{CaFG, derivation, JS3, JS, JS2, JS4}. We make another contribution but more in the line of the classic Kalton-Peck paper \cite{KaPe}. Indeed, the Kalton-Peck space $Z_2$ will be the center of gravity of our paper. This space may be simply defined as above as certain set of couples $(x,y)$ where $x,y\in \ell_2$, see \cite{KaPe}. The natural Hilbert copy is given by those vectors of the form $(x,0)$ and the quotient of $Z_2$ by this copy is again the Hilbert space $\ell_2$; the vectors $(0,x)$ span an Orlicz space usually denoted $\ell_M$. The so-called Kalton-Peck map $\mathcal K(x)=x\cdot\log|x|$ for $\|x\|=1$ (and extending by homogeneity) is responsible for the norm of $Z_2$: $$\|(x,y)\|_{Z_2}=\|x-\mathcal K(y)\|_{\ell_2}+\|y\|_{\ell_2}.$$ In the differential approach of Nigel, the space $Z_2$ appears as $Z(X)$ for $X=c_0$ or any $X=\ell_p$ for $1\leq p\neq 2<\infty$. This is, it appears as a derivation of any of these classic sequence spaces, except for the Hilbert copy. In this case, it is easy to check that $Z(\ell_2)=\ell_2$. We are very interested in the fact that Kalton and Peck developed a scheme of results for $Z_2$ showing that it is a very rigid space. For example:
\begin{itemize}
\item[($K_1$)] The natural quotient map $Q:Z_2\longrightarrow \ell_2$ is strictly singular.
\item[($K_2$)] Every normalized basic sequence in $Z_2$ has a subsequence equivalent to either $\ell_2$ or the Orlicz space $\ell_M$.
\end{itemize}

The choice of these claims is not by chance, we will see at the end of the section that both statements are actually connected.

Recently, Castillo, Ferenczi and Gonz\'alez dealt in \cite{CaFG} with $Z(X)$ where $X$ is a reflexive asymptotic $\ell_p$-space for $p\neq 2$. They showed that the natural quotient map of $Z(X)$ onto $\ell_2$ is strictly singular, meaning that the quotient map is never an isomorphism when restricted to a Hilbert copy of $Z(X)$, see \cite[Proposition 6.3.]{CaFG}.  The twisted Hilbert $Z(\mathcal T)$, where $\mathcal T$ denotes the Tsirelson space, is a good example of this phenomena as it is remarked in \cite{CaFG}. Recall that $\mathcal T$ is well known to be an asymptotic $\ell_1$-space; this is, for every $n\in \mathbb N$, every $n$ normalized blocks of $\mathcal T$ starting down far enough are equivalent to the basis of $\ell_1^n$. To sum up, the authors of \cite{CaFG} were interested in $(K_1)$ of the Kalton-Peck figurehead for $p\neq 2$. They left open the cases of $c_0$ and $p=2$ and no other claim concerning the isomorphic structure of such twisted Hilbert spaces is given in \cite{CaFG}. In particular, there are no results concerning $(K_2)$. 

 In \cite{JS}, the author showed that $Z(\mathcal T^2)$ is a weak Hilbert space which is an important class of spaces introduced by Pisier in \cite{Pi}. This twisted Hilbert space corresponds to the case $p=2$ since $\mathcal T^2$ is an asymptotic $\ell_2$-space. The isomorphic structure of this space is better understood as it is basically incomparable with the Kalton-Peck space or the Enflo-Lindenstrauss-Pisier example. However, facts are that it is not known how rigid the space $Z(\mathcal T^2)$ could be. The interested reader may check \cite{JS3} for more new examples of twisted Hilbert spaces.
 
 The purpose of this paper is to describe the rigidity for a twisted Hilbert space $Z(X)$ induced by an asymptotic $c_0$ or $\ell_p$-space $X$ for $1\leq p<\infty$, thus extending the Kalton-Peck results. For $p=2$ we show, for example, that $Z(\mathcal T^2)$ has a natural strictly singular quotient map onto $\ell_2$ which fills the gap left by Castillo \textit{et al}. in \cite{CaFG}. In general, we shall answer the following questions for such $X$:
 \begin{enumerate}
 \item[(A)] Is the natural quotient map from $Z(X)$ onto $\ell_2$ strictly singular?
 \item[(B)] Characterize the (normalized) basic or weakly null sequences in $Z(X)$.
 \end{enumerate}
A key step to answer question (A) is to isolate a rule we term as the Principle of Small Perturbations for a twisted Hilbert space $Z$: \textit{``All the Hilbert copies in $Z$ are basically a small perturbation of vectors of the form $(x,0)$."}

To answer $(B)$ we deal with spreading models. Spreading models were introduced by Brunel and Sucheston in the early 70's. Once the idea burst into the mathematical scene, it established at the heart of Banach space theory. The fundamental discovery of Brunel and Sucheston was that every bounded sequence $(x_j)$ in a separable Banach space $X$ generates a norm which is given asymptotically. 

In particular, if $X$ denotes an asymptotic $\ell_p$-space for $p\neq 2$, we have the following answer to (B): \textit{The only spreading models of $Z(X)$ are $\ell_2$ and the Orlicz space $\ell_M$.} We shall also introduce the notion of twisted spreading model showing that the above can be written as ``\textit{The only twisted spreading model of $Z(X)$ is $Z_2$}."

We will see in the last section that, as in the previous claim (in its twisted form), if it happens that some twisted Hilbert space $Z$ has only $Z_2$ as spreading model, then such $Z$ must satisfy the Principle of Small Perturbations. We recall that such rule is in turn responsible for the quotient map in $(A)$ to be singular. So questions $(A)$ and $(B)$ are actually deeply related. Indeed, it is the notion of maximal centralizer (see Section \ref{sectionmM} for definitions) which is truly connecting both statements. This notion is our starting point.

\subsection{Roadmap of the paper.}

The paper contains three parts.

\subsubsection*{\textbf{Part I}}\textit{ The Kalton-Peck circle of ideas for asymptotic $\ell_p$-spaces: The basics.}
\begin{itemize}
\item[(S3)] We introduce the notion of maximal centralizer which corresponds to $p\neq 2$ and mimics the Kalton-Peck centralizer asymptotically. This is the key for (S10).
\item[(S4)] We state The Principle of Small Perturbations.
\end{itemize}
\subsubsection*{\textbf{Part II}} \textit{Spreading models of a twisted Hilbert space with two good examples.}
\begin{itemize}
\item[(S5)] We introduce spreading models for a twisted Hilbert space.
\item[(S6)] We discuss the spreading models of the example $Z(\mathcal T^2)$ that corresponds to the case $p=2$.
\item[(S7)] The Principle of Small Perturbations for $Z(\mathcal T^2)$ appears using a tool from the previous section.
\item[(S8)] We discuss the spreading models of the example $Z(\mathcal T)$ that shows us how works the case $p\neq 2$ in a concrete situation.
\end{itemize}
\subsubsection*{\textbf{Part III}} \textit{The Kalton-Peck theory for asymptotic $\ell_p$-spaces $X$.}
\begin{itemize}
\item[(S9)] We develope the spreading models of $Z(X)$ for $p\neq 2$ which is only the abstract version of (S8).
\item[(S10)] We prove a sort of converse to the results of (S9) and relate these spreading models with the Principle of Small Perturbations and with the maximal centralizers of (S3) closing our circle of ideas.
\end{itemize}

\section{Background}
 In this paper we shall use the terminology commonly used in Banach space theory as it appears in the book of Albiac and Kalton \cite{AK}. The word ``space" means ``infinite dimensional Banach space" unless specified otherwise. We will denote by $[x_j]_{j\in A}$ the closed linear span of a sequence $(x_j)_{j\in A}$ in a space $X$. If $X$ has an unconditional basis, and no confusion arises, $(e_j)_{j=1}^{\infty}$ will be such basis. In this case, we write as usual $\supp(x)$ for the support of $x\in X$.  And, as usual too, given two spaces $X,Y$, we write $X\approx Y$ if they are isomorphic.
 
 In many proofs, $\eta$ and $C$ denote ``generic" constants, which need not have
the same value throughout the proof. In chains of inequalities, we will often use
$C,C_1,C_2,...$ to avoid confusion. Also, given two sequences of real numbers $(a_j)_{j=1}^{\infty}$ and $(b_j)_{j=1}^{\infty}$, we write $a_j \sim b_j$  if there is $C>0$ such that $C^{-1}\cdot b_j\leq a_j\leq C\cdot b_j$ for all $j\in \mathbb N$.
 

The main examples of this paper are based on Tsirelson's space $\mathcal T$ for which the book of Casazza and Shura \cite{CS} provides a comprehensive study. To introduce Tsirelson's space we need to give inductively a sequence of norms. Pick an element $x$ in $c_{00}$, the vector space of all finitely supported sequences, and define:
$$\|x\|_0=\max |x_j|,$$
$$\|x\|_{n+1}= \max \left \{ \|x\|_n, \frac{1}{2} \max \left [ \sum_{j=1}^k \|A_j x\|_n   \right ] \right\} , \;\;(n\geq 0),$$
where the inner max is taken over all choices of finite subsets as $k$ varies and such that $k\leq A_1<...<A_k$. We have written $A_j<A_{j+1}$ as usual for $\max A_j <\min A_{j+1}$ and $A_jx$ for the natural restriction of $x$ to the coordinates of $A_j$. Tsirelson's space $\mathcal T$ is the completion of $c_{00}$ under the norm $\|x\|:=\lim_{n\to \infty}\|x\|_n$. We will denote $(t_j)_{j=1}^{\infty}$ and $(t_j^*)_{j=1}^{\infty}$ the natural basis in $\mathcal T$ and $\mathcal T^*$ respectively. The $p$-convexification of $\mathcal T$ is denoted $\mathcal T^p$ and is the completion of $c_{00}$ under the norm
$$\left \|\sum x_j t_j \right \|_{\mathcal T^p}:=\left \| \sum |x_j|^p t_j \right \|_{\mathcal T}^{\frac{1}{p}},$$ where $1\leq p<\infty$.
In the particular case $p=2$, the space $\mathcal T^2$ is a weak Hilbert space which is an important and special class of spaces introduced by Pisier in \cite{Pi}. 
\begin{definition}\label{weakhilbert}
We say $X$ is a weak Hilbert space if there is $0<\eta<1$ and a constant $C$ with the following property: every finite-dimensional subspace $E$ of $X$ contains a subspace $F\subseteq E$ with $\dim F\geq \eta\cdot \dim E$ such that $d(F,\ell_2^{\dim F})\leq C$ and there is a projection $P:X\to F$ with $\|P\|\leq C$, where $d$ denotes the Banach-Mazur distance. 
\end{definition}
The definition above is not the original one but is chosen out among the many equivalent characterizations given by Pisier \cite[Theorem 12.2.(iii)]{Pi}. Since the property ``to be a weak Hilbert space" passes to subspaces, quotients and duals (\cite{Pi}), we find that $(\mathcal T^2)^*$ is also a weak Hilbert space.
There is a couple of facts on the structure of $\mathcal T^p$ and its dual that will be crucial for us:
\begin{itemize}
\item The spaces $\mathcal T^p$ are \textbf{asymptotic $\ell_p$-spaces}.
\end{itemize}
 This means that there is $C>0$ so that for any $n\in \mathbb N$ and any normalized blocks $n<u_1<...<u_n$, we find that $(u_j)_{j=1}^n$ is $C$-equivalent to the unit vector basis of $\ell_p^n$. Roughly, every $n$ consecutive blocks starting down far enough are just $\ell_p^n$. The same definition works for a space with basis.
\begin{itemize}
\item The space $\mathcal T^*$ is an \textbf{asymptotic $c_0$-space}.
\end{itemize}
We shall only need that for any $n\in \mathbb N$, any normalized blocks $n<u_1^*<...<u_n^*$ and scalars $\alpha_1,...,\alpha_n$ with $|\alpha_j|\leq 1$ for $j\leq n$ we have that 
\begin{equation}\label{formulilla}
\left\|\sum_{j=1}^n \alpha_ju_j^*\right\|_{\mathcal T^*}\leq 2.
\end{equation}
This follows easily from a standard duality argument. Indeed, pick $x\in \mathcal T$ with $\|x\|=1$ norming (\ref{formulilla}). There is no loss of generality assuming, using the unconditionality of the basis of $\mathcal T$, that $\supp x \subseteq \cup_{j=1}^n \supp(u_j^*)$. But it happens that if $A_j:=\supp(u_j^*)$, then $(A_j)_{j=1}^n$ is admissible, and thus
\begin{eqnarray*}
\left\| \sum_{j=1}^n \alpha_ju_j^* \right\|_{\mathcal T^*}&=&\left|\sum_{j=1}^n \alpha_ju_j^*(x)\right|\\
&=&\left|\sum_{j=1}^n \alpha_ju_j^*(A_j(x))\right|
\leq \sum_{j=1}^n \|A_j(x)\|_{\mathcal T}\leq 2\|x\|_{\mathcal T}\leq 2,
\end{eqnarray*}
where the previous to the last inequality holds by the very definition of the norm in $\mathcal T$. A very similar duality argument shows that $(\mathcal T^2)^*$ is also an asymptotic $\ell_2$-space.

Many of the claims of this paper are stated in the language of spreading models which fits perfectly when handling with $Z(\mathcal T)$ and $Z(\mathcal T^2)$. This language was introduced by Brunel and Sucheston as an application of Ramsey’s theorem. The fundamental discovery of Brunel and Sucheston was that every bounded sequence $(x_j)$ in a separable Banach space $X$ admits a subsequence $(x'_{j})$ such that given $n\in \mathbb N$ and scalars $a_1,...,a_n$ there is a number $L(a_1,...,a_n)$ verifying that for every $\varepsilon>0$, we may find $\nu\in \mathbb N$ with $$\left | \left\|\sum_{j=1}^na_jx'_{k_j}\right\|-L(a_1,...,a_n)\right|<\varepsilon$$ for all $\nu \leq k_1\leq ...\leq k_n$. Therefore, for each $a_1,...,a_n$ the formula
$$\left\|a_1e_1+...+a_ne_n\right\|:=L(a_1,...,a_n)$$ defines a semi-norm. If  $(x_j)$ admits no converging subsequences, then it is easy to show that the semi-norm is actually a norm \cite[Proposition 2, page 8]{BeLa}. The completion of $c_{00}$ with the norm described above is a Banach space called a \textbf{spreading model} of $X$ generated by the sequence $(x_j)$. And the sequence $(e_j)_{j=1}^{\infty}$ is termed as the \textbf{fundamental sequence} of the spreading model. If moreover, the sequence $(x_j)$ is weakly null and normalized then the fundamental sequence $(e_j)_{j=1}^{\infty}$ is $1$-supression unconditional \cite[Proposition 1, page 24]{BeLa}.

 
 For the topic of spreading models, we shall follow mainly the exposition of Beauzamy and Laprest\'e in \cite{BeLa}. 

\subsection{How to start our own twisted Hilbert space}\label{howto}We make use of a nowadays standard method due to Kalton \cite{ka} to construct twisted Hilbert spaces. It is only required a space $X$ with an unconditional basis. The key of the process lies on the idea of centralizer which is certain nonlinear map associated to such $X$. The final target is to get such map and thus the norm on the twisted Hilbert space. We would like to recall here some basic points.

 A homogeneous map $\Omega:\ell_2\longrightarrow \omega$ is called a \textbf{centralizer} if there is a constant $C>0$ so that: 
\begin{equation}\label{centralizer}
\| \Omega(ax)-a\Omega(x)\|_{\ell_2}\leq C\|a\|_{\infty}\|x\|_{\ell_2},\;\;\:\;a\in \ell_{\infty},x\in \ell_2,
\end{equation}
where $\omega$ denotes the vector space of all complex scalar sequences with the topology of pointwise convergence. The reason why $\Omega$ produces a norm or is related to a twisted Hilbert space is not clear at this stage and will take us some steps.

\textit{1. How do we obtain a centralizer starting from $X$?}

 A standard way is through complex interpolation for which we follow the approach of Kalton and Montgomery-Smith \cite{KM} while the classic reference is the book of Bergh and L\"ofstr\"om \cite{BeLo}. Let us fix the basics on interpolation here. Let $X,Y$ be spaces with a joint unconditional basis and natural inclusions into $\omega$; so we are identifying $X,Y$ as sequence spaces. Interpolation deals with the vector space $\mathcal F_{\infty}(X,Y)$ of all functions $F:\mathbb S\to \omega$, that are bounded and continuous on the strip $$\mathbb S=\{ z: 0\leq Rez \leq 1 \},$$
and analytic on the open strip $\{z: 0< Re z <1 \},$ and moreover, the functions $t\in \mathbb R\to F(it)\in X$ and $t\in \mathbb R\to F(1+it)\in Y$ are bounded and continuous functions. Such vector space is a Banach space when is endowed with the norm $$\|F\|_{\mathcal F_{\infty}}=\max \left ( \sup_{t\in \mathbb R} \| F(it)\|_{X},   \sup_{t\in \mathbb R} \| F(1+it)\|_{Y}\right).$$
 The definition may look like quite technical and a bit artificial. It is not much important here since the core for us is the interpolation space $(X,X^*)_{1/2}$ which has a norm equivalent to $\ell_2$ under our assumptions, see e.g. \cite{CoS} and \cite{Wat}. This isomorphism connects our choice $X$ with the Hilbert copy $\ell_2$ for the first time.  The space $(X,X^*)_{1/2}$ consists of all $x\in \omega$ such that $x=F(1/2)$ for some $F\in \mathcal F_{\infty}(X,X^*)$; so we may see $\mathcal F_{\infty}$ as a kind of structure over $\ell_2$. Indeed, we have a natural quotient map $\delta_{1/2}:\mathcal F_{\infty}(X,X^*) \to  \ell_2$ given by $\delta_{1/2}(F)=F(1/2)$ and we may think in $\ell_2$ as endowed with the quotient norm $$\inf \{ \|F\|_{\mathcal F_{\infty}} : F(1/2)=x\}.$$

\textit{2. Where is the centralizer?}

 Once $X$ and $\ell_2$ are connected, we must ``derive" the interpolation formula $(X,X^*)_{1/2}=\ell_2$ which is to derive this connection. To do so, we must fix a constant $C>1$ and for each $x\in \ell_2$ pick a map $B(x)\in \mathcal F_{\infty}(X,X^*)$ with $B(x)(1/2)=x$ and $\|B(x)\|_{\mathcal F_{\infty}}\leq C\|x\|_{\ell_2}$. There is no loss of generality assuming that our ``$B$" is an homogeneous map. A centralizer $\Omega:\ell_2\longrightarrow \omega$ comes defined as $$\Omega(x)=\delta'_{1/2}B(x).$$ 

\textit{3. We have a centralizer. What is next?}

 Centralizers and twisted Hilbert spaces are naturally related. A compact way to describe twisted Hilbert spaces is using a short exact sequence which is a diagram like
$$
\begin{CD} 0 @>>>\ell_2@>j>>Z(X)@>Q>>\ell_2@>>>0,
\end{CD}
$$ where $j(x)=(x,0)$ and $Q(x,y)=y$. The formal definition claims that the morphisms are linear and continuous and such that the image of each arrow is the kernel of the next one. This last condition compacts the whole information: that $\ell_2$ is a subspace of $Z(X)$ through $j$ and thanks to
the open mapping theorem we find that $\ell_2$ is isomorphic to the quotient $Z(X)/j(\ell_2)$. We thus label $Z(X)$ as a \textbf{twisted Hilbert space}. 

\textit{4. Setting up our $Z(X)$.}

 Once we have $\Omega$, pick the set of couples $(x,y)\in \omega\times \omega$ for which the quasi-norm \begin{equation}\label{quasinorm}
\|(x,y)\|:=\|x-\Omega(y)\|_{\ell_2}+\|y\|_{\ell_2},
\end{equation}
is finite. The easy part is that the quasi-norm produces a short exact sequence as before \cite[Page 1159]{KM}; indeed, pick $j(x)=(x,0)$ and $Q(x,y)=y$. The non-easy part is that the quasi-norm is equivalent to a norm by a result of Kalton \cite{ka3}. The interested reader in the exact but theoretical shape of such norm is encouraged to read \cite[Proposition 7.2.]{CCS}. Anyway, once both details are checked we are done.

\textit{5. Our twisted Hilbert space $Z(X)$ is well defined.} 

It will be useful to recall that a centralizer $\Omega$ is \textbf{bounded} if $\Omega(y)\in \ell_2$  for every $y\in \ell_2$ and there is $C>0$ such that $\|\Omega(y)\|_{\ell_2}\leq C \|y\|_{\ell_2}$ in which case (\ref{quasinorm}) is equivalent to the product norm. It is clear that there is some ambiguity in the choice of $\Omega$ but the difference of any two choices is a bounded centralizer what implies that the quasi-norms (\ref{quasinorm}) (and thus the equivalent norms) are equivalent.

\textit{6. Checking out our $Z(X)$.}

 Once $Z(X)$ has been constructed it is a legitimate question whether $Z(X)$ is again Hilbert or not. The answer to this is given by Kalton's uniqueness theorem \cite[Theorem 7.6.]{ka}. In its simplest form works as:
$$Z(X)\approx \ell_2\Longleftrightarrow X\approx \ell_2.$$
We give some hints to the interested reader about how to deduce the claim from \cite[Theorem 7.6.]{ka}. It is easy to show that if $Z(X)\approx \ell_2$, then the corresponding centralizer $\Omega$ is bounded, see e.g. \cite{JS2} for a detailed argument. If we pick the new couple $Y_0=Y_1=\ell_2$ whose centralizer is $0$, we are saying that $\Omega$ and $0$ are equivalent in the language of \cite{ka} (see also the next section). Since by \cite[Theorem 7.6.]{ka}, the spaces representing a centralizer are uniquely determined, we find $X\approx Y_0=\ell_2$.

\textit{7. An example to keep in mind.}

  The canonical example where all these ideas crystallize is the Kalton-Peck space which corresponds to $Z(c_0)=Z_2$. It is well known that the centralizer comes defined for a norm one vector $x=\sum_{j=1}^{\infty}x_je_j$ as 
\begin{equation}\label{kpmap}
\mathcal K(x)=\sum_{j=1}^{\infty}x_j\log|x_j|e_j,
\end{equation} with the agreement that $0\cdot \log 0=0$. The map (\ref{kpmap}) is the so-called Kalton-Peck centralizer \cite{KaPe}. That the space $Z_2$ contains the following sequence space with symmetric basis
\begin{equation*}\label{orlicz} 
\ell_{M}:=[(0,e_j)]_{j=1}^{\infty}
\end{equation*}
will be a crucial point for us in many arguments. Furthermore, the space $\ell_M$ is an Orlicz space \cite{KaPe}. 

 One last important fact for us about $Z(X)$ is the following.
\begin{proposition}\label{lambda} If $(e_j)_{j=1}^{\infty}$ denotes the unconditional basis of $\ell_2$, then
\begin{enumerate}
\item[$(i)$] The sequence $\{(e_j,0),(0,e_j)\}_{j=1}^{\infty}$ is a basis for $Z(X)$.
\item[$(ii)$] The sequence $\{(0,e_j)\}_{j=1}^{\infty}$ is unconditional.
\end{enumerate}
\end{proposition}
The first part is proved by adapting the proof of \cite[Theorem 4.10]{KaPe} while the second follows picking $a\in \{-1,1\}^{n}$ in (\ref{centralizer}).
\section{Minimal and Maximal centralizers}\label{sectionmM}
This section is rather technical and the main points under consideration are the idea of minimal/maximal centralizer and the nontrivial facts that the centralizers corresponding to $Z(\mathcal T^2)$ and $Z(\mathcal T)$ are minimal and maximal respectively. At the end of the paper, the role of maximal centralizer will be the protagonist showing that it is responsible for the stability of the extreme properties of $Z_2$ by ``twisted" spreading models. So, this notion is truly beating at the heart of $Z_2$. A little before, in the next section, it will play the role of the perfect example in The Principle of Small Perturbations which is a leading rule of this paper. The reader used to working with spreading models will nose out something familiar in the definition of minimal/maximal centralizer.

Gowers's favourite pedadogical principle is \textit{examples first} \cite{G}, so let us sketch the ideas behind our definitions. We deal first with the twisted Hilbert space $Z_2$ which is governed by the Kalton-Peck map $\mathcal K$. Observe that we have $\left\|\mathcal K\left(\sum_{j=1}^n e_j\right)\right\|_2=\sqrt{n}\log\sqrt{n}$ and $\mathcal K(e_j)=0$ for $j\leq n$. On the other hand, for the twisted Hilbert space $\ell_2$, which is completely different from $Z_2$, we may pick as centralizer $\Omega=0$ or any linear map $L$. The point is that we are allowed to perturbate the linear map $L$ by a bounded map, say $B$, and obtain a centralizer which is equivalent to the previous one. This amounts to say that the corresponding twisted Hilbert space is again isomorphic to $\ell_2$. Either way, if we try to quantify the nonlinearity of $\mathcal K$ and $\Omega=L+B$, a possible and natural way is:

$$\left\|\mathcal K\left(\sum_{j=1}^n e_j\right)-\sum_{j=1}^n\mathcal K(e_j)\right\|_2=\sqrt{n}\log\sqrt{n},$$
$$\left\|B\left(\sum_{j=1}^n e_j\right)-\sum_{j=1}^nB(e_j)\right\|_2\leq 2 \|B\|\sqrt{n}.$$
Thus, the factor $\log\sqrt{n}$ is reflecting the separation between the trivial case $\ell_2$ from the extreme case $Z_2$. In the case of $\mathcal K$ the same estimate as above holds for blocks as we are about to check, so our abstract and spreading version of the behaviour of $\mathcal K$ goes as follows.
\begin{definition}
We say that a centralizer $\Omega:\ell_2\to \omega$ is \textbf{maximal} if whenever we have a sequence of semi-normalized and disjoint blocks $(u_j)$ in $\ell_2$, there is a subsequence $(u_j')$ such that for every $n\in \mathbb N$ there is $\nu\in \mathbb N$ satisfying
$$\frac{\|\Omega(\sum_{j=1}^n u'_{k_j})-\sum_{j=1}^n \Omega(u'_{k_j})  \|}{\sqrt n}\longrightarrow \infty\;\;\textit{as}\;\;n\to \infty,$$
whenever $\nu\leq k_1<...<k_n$.
\end{definition}
%
In the opposite side of the spectrum, we have the abstract behaviour of a trivial centralizer in this setting.
\begin{definition}
We say that a centralizer $\Omega:\ell_2\to \omega$ is \textbf{minimal} if there is a constant $C>0$ such that whenever we have a sequence of $C_1$-bounded and disjoint blocks $(u_j)$ in $\ell_2$, we find a subsequence $(u_j')$ satisfying: for every $n\in \mathbb N$ there is $\nu\in \mathbb N$ such that
$$\left\|\Omega\left(\sum_{j=1}^n u'_{k_j}\right)-\sum_{j=1}^n \Omega(u'_{k_j})  \right\|\leq C\cdot C_1\sqrt n,$$
for every $\nu\leq k_1<...<k_n$.
\end{definition}

The notion of equality between two centralizers defines $\Omega$ and $\Omega'$ to be \textbf{equivalent} if the difference is a bounded centralizer; in particular $\Omega$ and $\Omega'$ produce equivalent quasi-norms and hence isomorphic twisted Hilbert spaces. Thus, the next lemma is no more than a simple exercise.
\begin{lemma} Assume that $\Omega$ is equivalent to $\Omega'$, then
\begin{enumerate}
\item The centralizer $\Omega$ is minimal if and only if $\Omega'$ is minimal.
\item The centralizer $\Omega$ is maximal if and only if $\Omega'$ is maximal.
\end{enumerate}
\end{lemma}
So, our definitions are stable through the natural equivalence relation. Following now Halmos \cite{Hal}, our first job after a definition is to build an example but we just need to confirm that $\mathcal K$ is maximal and $\Omega=0$ minimal.
\begin{example} The Kalton-Peck centralizer $\mathcal K$ is maximal.
\end{example}
\begin{proof}
Let $(u_j)$ be any sequence of semi-normalized blocks, so that for some $\alpha,\beta$ we have $\alpha\leq \|u_j\|\leq \beta$ for $j=1,2,...$. Then, it is proved in \cite[Theorem 6.5, page 28, line 17]{KaPe} that there is a constant $c(\alpha,\beta)$  so that 
$$c(\alpha,\beta)^{-1}\|(a,b)\|_{Z_2}\leq \|U(a,b)\|$$ for any $(a,b)\in c_{00}\times c_{00}$, where $$\|U(a,b)\|=\left\| \sum_{j=1}^{n} a_ju_{j}+\sum_{j=1}^n b_j\mathcal K(u_{j})-\mathcal K\left(\sum_{j=1}^n b_ju_{j}\right) \right\|_2+\left\| \sum_{j=1}^n b_ju_{j}\right\|_2.$$
Picking $a_1=...=a_n=0$ and $b_1=...=b_n=1$, we easily find using the previous estimates and the definition of $\mathcal K$ that $$\|U(0,1)\|\geq c(\alpha,\beta)^{-1}\sqrt{n}\log\sqrt{n},$$
and the claim follows trivially. We do not need to pass to a further subsequence.
\end{proof}
It is clear that $\Omega=0$ is minimal. However, there are many others minimal and maximal centralizers. The proof of this claim is not a trivial task since we do not know in general the shape of other centralizers. Therefore, we have no way to compute explicitly the exact value of the expressions involved in the definition of minimal/maximal  centralizer as we did for $\mathcal K$. So the proof of our claim requires some work.

 We fix a piece of notation for the rest of the section. Let $X$ be a space with an unconditional basis which is asymptotic $\ell_1$ for vectors with disjoint support; keep in mind that $X=\mathcal T$ is our model. For simplicity we shall assume that $X$ is reflexive while the Radon-Nikodym property is enough for most of the results; the interested reader may check e.g. \cite{BenLin} for more information on the Radon-Nikodym property. Let $X^p$ denote the $p$-convexification of $X$ and $\Omega_p$ will stand for a centralizer arising from the interpolation formula $(X^p,(X^p)^*)_{1/2}=\ell_2$ as described in Section \ref{howto} (see \textit{Where is my centralizer?}) and where $1\leq p<\infty$.  The following observation was communicated to us by Castillo \cite{Ca}:
\begin{lemma}\label{castillo}The centralizers
$$\Omega_p\;\;\textit{and}\;\; \frac{1}{p}\Omega_1+\left(1-\frac{1}{p}\right)\mathcal K$$
are equivalent for  $1\leq p < \infty$.
\end{lemma}
\begin{proof} Our proof is different from the one in \cite{Ca}. Our idea is based on the fact that centralizers $\Omega$ given by interpolation $X=(X_0,X_1)_{\theta}$ of a couple $(X_0,X_1)$ are uniquely determined, up to equivalence, by the entropy functions $$\phi_{X_k}(u)=\sup_{\|x\|_{X_k}\leq 1}\sum_{j=1}^{\infty}u_j\log|x_j|\;\;\textit{for}\;\; k=0,1,$$ see e.g. \cite{KM}. Here, determined means that there is $C>0$ so that for $x\in c_{00}$ and $x^*\in X^*$, we have $$\left|\langle\Omega(x),x^*\rangle-(\phi_{X_1}(xx^*)-\phi_{X_0}(xx^*))\right|\leq C\|x\|_X\|x^*\|_{X^*}.$$ We will write this claim in short as $\Omega=\phi_{X_1}-\phi_{X_0}$. Now, the formula $X^p=(X,c_0)_{1/p}$ is easy to establish using the Radon-Nikodym property of $X$. Indeed, we may invoke \cite[Theorem 4.6.]{KM} to find $(X,c_0)_{1/p}=X^{1/p}c_{0}^{1-1/p}$ isometrically from where the claim follows trivially. Since $X$ has the Radon-Nikodym property, we may indentify the dual of $(X,c_0)_{1/p}$ with $(X^*,\ell_1)_{1/p}$, see also \cite{KM}. Since the entropy function linearizes interpolation, this is $\phi_X=\theta\phi_{X_0}+(1-\theta)\phi_{X_1}$, we find 
\begin{eqnarray*}
\Omega_p&=& \phi_{(X^p)^*}-\phi_{X^p}\\
&=&\left[\frac{1}{p}\phi_{X^*}+\left(1-\frac{1}{p} \right)\phi_{\ell_1}\right]-\left[\frac{1}{p}\phi_{X}+\left(1-\frac{1}{p} \right)\phi_{c_0}\right]\\
&=&\frac{1}{p}\left[\phi_{X^*}-\phi_{X}\right]+\left(1-\frac{1}{p} \right)\left[\phi_{\ell_1}-\phi_{c_0}\right]\\
&=& \frac{1}{p}\Omega_1+\left(1-\frac{1}{p} \right) \mathcal K.
\end{eqnarray*}
The equality above claims truly that for $x\in c_{00}$ and any $x^*\in \ell_2$, we have
$$\left|\left\langle \Omega_p(x)-\left[ \frac{1}{p}\Omega_1(x)+\left( 1-\frac{1}{p}\right)\mathcal K(x)\right], x^*\right\rangle \right|\leq C\|x\|\|x^*\|,$$
for a suitable constant $C>0$. For a fixed $x\in c_{00}$, taking supremum over all $x^*$ of norm $1$ in the inner expression above we have
$$\left\|\Omega_p(x)-\left[ \frac{1}{p}\Omega_1(x)+\left( 1-\frac{1}{p}\right)\mathcal K(x)\right]\right\|\leq C\|x\|,$$
which is to say that the centralizers in the statement of the Lemma are equivalent.
\end{proof}

A simple computation gives us now
\begin{corollary}\label{p2} The centralizers 
$$\Omega_p\;\;\textit{and}\;\;\frac{2}{p}\Omega_2+\frac{p-2}{p}\mathcal K$$
are equivalent for $1\leq p <\infty$.
\end{corollary}
\begin{proof}
Let us write the claim of Lemma \ref{castillo} for $p=2$ and an arbitrary $p$ as
$$\Omega_2\cong\frac{1}{2}\Omega_1+\frac{1}{2}\mathcal K$$
and 
$$\Omega_p\cong\frac{1}{p}\Omega_1+\left(1-\frac{1}{p}\right)\mathcal K$$
respectively. If we clear $\Omega_1$ in the first expression and plug it into the second we are done.
\end{proof}
We shall infer now the existence of nontrivial maximal centralizers.
\begin{proposition}\label{hipotesis} Let $X$ be  a reflexive space with an unconditional basis. If $X$ is an asymptotic $\ell_1$-space for vectors with disjoint support, then we have
\begin{enumerate}
\item The centralizer $\Omega_2$ is minimal.
\item The centralizer $\Omega_p$ is maximal for every $p\neq 2$.
\end{enumerate}
\end{proposition}
\begin{proof}
The claim ``$\Omega_2$ is minimal" was proved by the author in \cite[Proposition 3]{JS2} for $\mathcal T^2$. The proof for $X^2$ is identical since the key of the argument in \cite{JS2} is based on the property $(H)$, which is, for every $C>1$ there is a constant $f(C)$ such that for any normalized $C$-unconditional basic sequence $u_1,...,u_n$ we have 
\begin{equation*}\label{Hache}
f(C)^{-1}\sqrt{n}\leq \left \|\sum_{j=1}^n u_j\right \|\leq f(C)\sqrt{n}.
\end{equation*}
It is well known that weak Hilbert spaces have the property $(H)$ 
\cite[Proposition 14.2.]{Pi}, so that $X^2$ and $(X^2)^*$ have property $(H)$ since they are weak Hilbert spaces. The reason for this is that $X^2$ is easily shown to be asymptotic $\ell_2$ for vectors with disjoint support and this implies to be weak Hilbert \cite[Theorem 3.1.]{ACK}. To prove the claim on $\Omega_2$, let $Y$ be for a moment a space with basis and let us write $u_1,...,u_n$ for $n$ blocks of such basis. We introduce the following parameter 
 $$D_n(Y)=\sup  \left \{    \|u_1+...+u_n\|_Y :u_1<...< u_n, \|u_j\|\leq 1, j\leq n \right\}.$$ 
 We may use the inequality \cite[Lemma 4.8.]{CaFG} with the parameter $D_n$ to find
 \begin{eqnarray*}
\left \| \Omega_2 \left(\sum_{j=1}^n u_j \right)- \sum_{j=1}^n\Omega_2 (u_j)\right\|&\leq& 6D_n(X^2)^{1/2}D_n((X^2)^*)^{1/2}\\
&+& \left \|  \log \frac{D_n(X^2)}{D_n((X^2)^*)} \sum_{j=1}^n u_j\right \|.
 \end{eqnarray*}
And by the property $(H)$, we have that $D_n(X^2)\sim \sqrt{n}\sim D_n((X^2)^*)$ so we are done.
 Using this, let us show the claim on $\Omega_p$. We just need to observe that $\mathcal K$ is maximal, $\Omega_2$ is minimal and use Corollary \ref{p2}. Indeed, we have that, up to the bounded factor of the equivalence, 
\begin{eqnarray*}
\left\|\Omega_p\left(\sum_{j=1}^n u_j\right)-\sum_{j=1}^n \Omega_p(u_j)\right\|&\geq& \frac{|p-2|}{p}\left\|\mathcal K\left(\sum_{j=1}^n u_j\right)-\sum_{j=1}^n \mathcal K(u_j)\right\|\\
&-&\frac{2}{p}\left\|\Omega_2\left(\sum_{j=1}^n u_j\right)-\sum_{j=1}^n \Omega_2(u_j)\right\|,
\end{eqnarray*}
from where the result follows for $p\neq 2$.
\end{proof}
\subsection{An application to the classification of twisted Hilbert spaces}
A straight application of maximal and minimal centralizers is to give a quick distinction of twisted Hilbert spaces. The following is the abstract version of the argument given in \cite[Corollary 1]{JS2} that generalizes in particular \cite[Proposition 4]{JS2}. 
\begin{proposition}\label{corcor}Let $Z_M,Z_m$ be twisted Hilbert spaces induced by a maximal centralizer $\Omega_M$ and a minimal centralizer $\Omega_m$, respectively. Then $Z_M$ is not isomorphic to a subspace of $Z_m$.
\end{proposition}
\begin{proof}
Let us assume that $T:Z_M\to Z_m$ is such embedding. Pick $\{(0,e_j)\}_{j=1}^{\infty}$ in $Z_M$ and using a standard gliding hump argument pass to a subsequence which we do not relabel so that $T(0,e_j)$ is a block basis $(a_j,b_j)$ in $Z_m$, where $\|a_j-\Omega_m(b_j)\|+\|b_j\|\leq 2\|T\|$ and the blocks $(a_j,b_j)$ are disjoint. Passing to further subsequences, there is no loss of generality assuming that $\{(0,e_j)\}$ and $\{(a_j,b_j)\}$ are the corresponding subsequences in the definition of maximal and minimal centralizer, respectively. Thus, given $n\in \mathbb N$, find the corresponding $\nu$ (for both sequences) and let $\nu<k_1<...<k_n$. We have
\begin{eqnarray*}
\left\|\left(0,\sum_{j=1}^n e_{k_j}\right)\right\|&\leq &C\left(\left\|\sum_{j=1}^n a_{k_j}-\Omega_m\left( \sum_{j=1}^n b_{k_j}\right)\right\|+ \left\| \sum_{j=1}^n b_{k_j} \right\|\right)\\
&\leq &C\left(\left\|\sum_{j=1}^n a_{k_j}-\sum_{j=1}^n \Omega_m\left( b_{k_j}\right)\right\|+ \left\| \sum_{j=1}^n b_{k_j} \right\|\right)\\
&+&C\left\|\sum_{j=1}^n \Omega_m\left( b_{k_j}\right)-\Omega_m\left( \sum_{j=1}^n b_{k_j}\right)\right\|\\
\end{eqnarray*}
Recalling that $T(0,e_j)=(a_j,b_j)$, we may still bound as
\begin{eqnarray*}
\left\|\left(0,\sum_{j=1}^n e_{k_j}\right)\right\|&\leq&4C\|T\|\sqrt{n}+C\left\|\sum_{j=1}^n \Omega_m\left( b_{k_j}\right)-\Omega_m\left( \sum_{j=1}^n b_{k_j}\right)\right\|\\
&\leq&C_1\sqrt{n},
\end{eqnarray*}
where the last inequality holds by the minimality of $\Omega_m$ since $\nu<k_1<...<k_n$. Therefore, normalizing the expression by $\sqrt{n}$ and using the maximality of $\Omega_M$, we reach a contradiction: $$C_1\geq \frac{\left\|\left(0,\sum_{j=1}^n e_{k_j}\right)\right\|}{\sqrt{n}}\geq \frac{\left\|\sum_{j=1}^n \Omega_M\left( e_{k_j}\right)-\Omega_M\left( \sum_{j=1}^n e_{k_j}\right)\right\|}{\sqrt n} \longrightarrow\infty,$$
as $n\to \infty$.
\end{proof}
We have now plenty of examples of non isomorphic twisted Hilbert spaces.
\begin{corollary}\label{corcor}Let $Z_M$ as above and let $X$ be a space with an unconditional basis that is an asymptotic $\ell_1$-space for vectors with disjoint support. Then $Z_M$ is not isomorphic to a subspace of $Z(X^2)$.
\end{corollary}
\begin{proof}
The claim ``$\Omega_2$ is minimal" is independent of the reflexivity of $X$.
\end{proof}
As an important particular case, let us isolate the following.
\begin{corollary}\label{maxmin} The twisted Hilbert space $Z(\mathcal T)$ is not isomorphic to a subspace of $Z(\mathcal T^2)$.
\end{corollary}
Clearly, the converse to Proposition \ref{corcor} does not hold since we always may take $Z_m:=Z(\ell_2)=\ell_2$ which embeds in any $Z_M$. In sharp constrast, we will develope quite a few tools to show in Proposition \ref{zt2zt} that $Z(\mathcal T^2)$ is not isomorphic to a subspace of $Z(\mathcal T)$.

\section{The Principle of Small Perturbations}\label{singaff}
The following rule is at the root under the many extreme properties of the Kalton-Peck space. But far from being an exclusive fact on $Z_2$, the rule holds true for many interesting twisted Hilbert spaces. It claims, roughly speaking, that if a twisted Hilbert space $Z$ satisfies it, then all the copies of $\ell_2$ in $Z$ are just a small perturbation of the natural copy of $\ell_2$ in the first coordinate.
\begin{psp} A twisted Hilbert space $Z$ satisfies it if for every sequence $(w_j)$ in $Z$ which is equivalent to the unit vector basis of $\ell_2$ there is a subsequence $(w_{k_j})$ of $(w_j)$ and blocks $(y_j,0)$ in $Z$ such that 
\begin{enumerate}
\item $\|(y_j,0)\|=1$ for every $j\in \mathbb N$.
\item $\|w_{k_j}-(y_j,0)  \|\leq 2^{-j}$ for every $j\in \mathbb N$.
\end{enumerate}
\end{psp}
As the perfect example, we have those twisted Hilbert spaces governed by a maximal centralizer.
\begin{proposition}\label{PSP}
Let $X$ be a space with an unconditional basis that is a reflexive asymptotic $\ell_1$-space for vectors with disjoint support. Then the twisted Hilbert space $Z(X^p)$ satisfies The Principle of Small Perturbations for $p\neq 2$.
\end{proposition}
\begin{proof} The argument is contained in the proof of \cite[Theorem 16.16.]{BenLin} of Benyamini and Lindenstrauss for $Z_2$. So, let $(w_j)$ be a sequence as above and, since $w_j\to 0$ weakly, pass to a subsequence (we do not relabel) so that there exists blocks $(x_j,u_j)$ with $$\|w_j-(x_j,u_j)\|\leq 2^{-j},\;\;j\in \mathbb N.$$ We claim that $\|u_j\|\to 0$. Let us observe that we always have the inequality
\begin{eqnarray*}
\left\| \left(\sum_{j=1}^n x_j,\sum_{j=1}^n u_j\right) \right\| &\geq& \left\| \sum_{j=1}^n x_j-\Omega_p\left( \sum_{j=1}^n u_j \right) \right\|\\
&\geq &  \left\| \Omega_p\left(\sum_{j=1}^n u_j\right)-\sum_{j=1}^n\Omega_p\left( u_j \right) \right\|-  \left\| \sum_{j=1}^n (x_j-\Omega_p\left(u_j\right) \right\|,
\end{eqnarray*}
which turns out to be impossible if $\|u_j\|\geq \eta >0$. Indeed, $$\left\| \sum_{j=1}^n \left(x_j-\Omega_p(u_j)\right)\right\|^2\leq \sum_{j=1}^n \|x_j-\Omega_p(u_j) \|^2\leq n,$$ while $$\left\| \left(\sum_{j=1}^n x_j,\sum_{j=1}^n u_j\right) \right\| \sim \sqrt{n}$$ because the blocks $(x_j,u_j)$ are a small perturbation of $(w_j)$. If we assume that $\|u_j\|\geq \eta >0$, then the blocks $u_j$ are semi-normalized and, by maximality, the inequality breaks down for $n$ large enough: divide by $\sqrt{n}$ and let $n\to \infty$. Thus, it must be $\|u_j\|\to 0$. So let us assume $\|u_j\|_2\leq 2^{-j}$, and hence 
\begin{eqnarray*}
\|w_j-(x_j-\Omega_p(u_j),0)\|&\leq& \|w_j-(x_j,u_j)\|+\|(\Omega_p(u_j),u_j)\|\\
&\leq& 2^{-j}+\|u_j\|_2\\
&\leq& 2^{-j+1}.
\end{eqnarray*}
Thus, pick $y_j:=x_j-\Omega(u_j)$. Replacing $w_j$ by $\lambda_j w_j$, where the sequence $\lambda_j$ satisfies $\inf_j \lambda_j>0$ and $\sup_j \lambda_j<\infty$, we may assume that $\|y_j\|_2=1$.
\end{proof}
The following result on strictly singular operators was already obtained by Castillo \textit{et al.} \cite[Proposition 6.3.]{CaFG}. Recall than an operator $T:X\to Y$ is \textbf{strictly singular} if no restriction to a infinite-dimensional closed subspace is an isomorphism.
\begin{corollary}\label{one} Let $X$ be a space with an unconditional basis that is a reflexive asymptotic $\ell_1$-space for vectors with disjoint support. Then, the natural quotient map of $Z(X^p)$ onto $\ell_2$ is strictly singular for every $p\neq 2$.
\end{corollary}
\begin{proof}
Let us recall first that every infinite-dimensional subspace $W$ of a twisted Hilbert space contains an isomorphic copy of $\ell_2$. Indeed, if $W\cap \Ker Q$ is infinite-dimensional, where $Q(a,b)=b$ denotes the natural quotient map onto $\ell_2$, this is clear. Otherwise, $Q$ is an isomorphism on a finite codimensional subspace of $W$ and we are done again.

 So let us assume, on the opposite, that there is a sequence $(w_j)$ equivalent to the unit vector basis of $\ell_2$ such that $$1=\|w_j\| \sim \|Q(w_j)\|.$$ Since, by The Principle of Small Perturbations, we always have, passing to a subsequence of $(w_j)$ if necessary, that $$\|Q(w_j)\|=\|Q(w_j)-Q(y_j,0)\|\leq \|Q\|2^{-j},$$
for certain normalized blocks $(y_j,0)$, we reach a contradiction.
\end{proof}
\section{Spreading models for twisted Hilbert spaces}\label{spreadmodel}
 The idea of introducing asymptotic $\ell_p$-spaces to produce twisted Hilbert spaces, as it was done in \cite{CaFG,JS}, was demanding in our opinion a study of spreading models in twisted Hilbert spaces. It is clear that every twisted Hilbert space admits $\ell_2$ as a spreading model just by picking the sequence $\{(e_j,0)\}_{j=1}^{\infty}$. However, this does not give any new information. On the other hand, we may pick $\{(0,e_j)\}_{j=1}^{\infty}$ that generates the Orlicz space $\ell_M$ in the Kalton-Peck space. The spreading model is again the Orlicz space  $\ell_M$ which is not a twisted Hilbert space. It is perhaps not surprising that $Z_2$ has a natural basis $\{(e_j,0),(0,e_j)\}_{j=1}^{\infty}$, so if split the basis of $Z_2$ into the odd and the even terms respectively, we get only the aforementioned, say $1$-dimensional spaces. The magic of $Z_2$ begins when we consider all the basis together. In the same line, the interesting part of $Z_2$ as a ``twisted" spreading model begins when we consider couples of vectors of the basis. This is, when we see $Z_2$ asymptotically in dimension $2$. Then, the only twisted spreading model of $Z_2$ is again $Z_2$, see Example \ref{example}. This resembles our idea that $Z_2$ is a very rigid space. In this sense, we shall deal with sequences of the form $$\{ (u_j,0),(\Omega(u_j),u_j)  \}_{j=1}^{\infty},$$ where $\Omega$ is the centralizer corresponding to our twisted Hilbert space and $(u_j)$ will denote for the rest of the section a sequence of normalized and consecutive disjoint blocks in $\ell_2$. We shall refer formally in Definition \ref{definition} to the asymptotic model it generates as a ``twisted spreading model".

  For the philosophical reader, twisted spreading models are a generalization of spreading models as, for example, the one of Halbeisen and Odell \cite{HalOdell}. They extend standard spreading models in a similar vein as we may consider $Z_2$ an extension of $\ell_2$. For the pragmatic reader, twisted spreading models (or spreading models of twisted Hilbert spaces) may be used to distinguish twisted Hilbert spaces or even subspaces of such twisted Hilbert spaces.

We start the engine as there is no other way with the twisted version of the Brunel-Sucheston procedure for extracting good subsequences (\cite{BrSu}) which relies again on Ramsey's theorem.  
\begin{proposition}\label{BrSu} Let $Z$ be a twisted Hilbert space with centralizer $\Omega$. Every sequence $\{ (u_j,0),(\Omega(u_j),u_j)  \}_{j=1}^{\infty}$ contains a subsequence $$\{ (u'_{j},0),(\Omega(u'_{j}),u'_{j})  \}_{j=1}^{\infty}$$ with the following property:\\
 For each $(a,b)\in c_{00}\times c_{00}$ there is a number $L(a,b)$ such that for every $\varepsilon>0$, there exists $\nu\in \mathbb N$ with $$\left | \left\|\sum_{j=1}^{\infty}a_j(u'_{n_j},0)+\sum_{j=1}^{\infty} b_j(\Omega(u'_{n_j}),u'_{n_j})  \right\|-L(a,b)\right|<\varepsilon$$ for all $\nu \leq n_1\leq n_2\leq...$.
\end{proposition}
\begin{proof}The process is exactly as the one given by Brunel-Sucheston \cite[Proposition 1]{BrSu} but just for couples.

\end{proof}
We define a (semi) quasi-norm on $c_{00}\times c_{00}$ as
\begin{eqnarray*}\label{qn}
 \left\|(a_1,...,a_n,b_1,...,b_n)\right\|&:=&\lims \left\|\sum_{j=1}^na_j(u'_{k_j},0)+\sum_{j=1}^n b_j(\Omega(u'_{k_j}),u'_{k_j})  \right\|\\
  &=& \lims \nabla(a,b)+\left\| \sum_{j=1}^n b_ju'_{k_j}\right\|_2,
\end{eqnarray*}
where $$\nabla(a,b)=\left\| \sum_{j=1}^n a_ju'_{k_j}+\sum_{j=1}^n b_j\Omega(u'_{k_j})-\Omega \left(\sum_{j=1}^n b_ju'_{k_j}\right) \right\|_2.$$
Let us recall that given $(a,b)\in Z$ and $(x,y)\in Z^*$, the natural duality pairing is given as $$\langle(x,y),(a,b)\rangle=\langle x,b \rangle+ \langle y, a\rangle,$$
see e.g. \cite{JS2}. Therefore, it is easy to see that the sequence $$\{(u'_{j},0), (\Omega(u'_{j}),u'_{j}) \}_{j=1}^{\infty}$$ is weakly null. Since it is normalized, it is not convergent being $0$ the only possible candidate to a limit. Thus, we may apply the standard argument of \cite[Lemma 2]{BrSu} also contained in \cite[Proposition 2, page 8]{BeLa}, to find that $$\|(a,b)\|=0 \Longrightarrow a=b=0.$$
\begin{definition}\label{definition}
We shall refer to the completion of $c_{00}\times c_{00}$ with the quasi-norm above as a \textbf{twisted spreading model} of $Z$ generated by $\{ (u'_{j},0),(\Omega(u'_{j}),u'_{j})  \}_{j=1}^{\infty}$  with fundamental sequence $\{(e_j,0),(0,e_j)\}_{j=1}^{\infty}$.
\end{definition}
%
%

The Kalton-Peck space is the test space for all the experiments in the twisted Hilbert space setting. It is a simple exercise using the shape of the Kalton-Peck map that the twisted spreading model of $\{(e_j,0),(0,e_j)\}_{j=1}^{\infty}$ in $Z_2$ is isomorphic to $Z_2$. It takes a calculation to show that indeed this is the only possibility.
\begin{example}\label{example} All the twisted spreading models of $Z_2$ are isomorphic to $Z_2$.
\end{example}
\begin{proof} Let $(u_j)$ be any sequence of semi-normalized blocks, so that for some $\alpha,\beta$ we have $\alpha\leq \|u_j\|\leq \beta$ for $j=1,2,...$. Then, it is shown in the proof of \cite[Theorem 6.5, page 28, line 17]{KaPe} that there are constants $c(\alpha,\beta)$ and $C(\alpha,\beta)$ so that 
$$c(\alpha,\beta)^{-1}\|(a,b)\|_{Z_2}\leq \|U(a,b)\| \leq C(\alpha,\beta)\|(a,b)\|_{Z_2},$$
for any $(a,b)\in c_{00}\times c_{00}$, where $$\|U(a,b)\|=\left\| \sum_{j=1}^{n} a_ju_{j}+\sum_{j=1}^n b_j\mathcal K(u_{j})-\mathcal K\left(\sum_{j=1}^n b_ju_{j}\right) \right\|_2+\left\| \sum_{j=1}^n b_ju_{j}\right\|_2.$$
Therefore, passing to the spreading subsequence of Proposition \ref{BrSu} gives again the Kalton-Peck norm using the same estimate as before. This is, the Kalton-Peck is invariant under twisted spreadings. 
\end{proof}
The claim of the example above fits with our idea that $Z_2$ is a very rigid space; a claim which is not surprising for the reader familiar with $Z_2$. Thus, asymptotically, the space $Z_2$ is rigid in dimension $2$ but not in dimension 1 in which we have two spreading models:
\begin{example}\label{example2} The spreading models of $Z_2$ are the Hilbert copy $\ell_2$ and the Orlicz space $\ell_M$.
\end{example}
\begin{proof}
Kalton and Peck proved that there are only two types of (weakly null) sequences in $Z_2$, the ones spanning a Hilbert copy and the ones spanning the Orlicz space $\ell_M$. Given any other sequence in $Z_2$, we may pass to a subsequence which is one as described before (\cite[Lemma 5.3.]{KaPe}). This is clearly enough.
\end{proof}
 This uniqueness of twisted spreading models is not an exclusive fact of the Kalton-Peck space since it also holds, for example, for $Z(\mathcal T^2)$ and $Z(\mathcal T)$ as we will see in the next sections. However, we will find a crucial difference between these two examples.
\section{The only spreading model of $Z(\mathcal T^2)$}\label{smzt2} 
We would like to recall that every subspace of a twisted Hilbert space contains an isomorphic copy of $\ell_2$ as we have seen in the proof of Corollary \ref{one}. This fact can be convoluted as every subspace of a twisted Hilbert space contains a further subspace whose spreading model is $\ell_2$. But, it is not true that every subspace contains a further subspace whose twisted spreading model is $\ell_2$. This is the case of $Z_2$ in which no twisted spreading model is isomorphic to $\ell_2$ as it was shown in Example \ref{example}. Once this point has been clarified, our aim is to show that in sharp contrast to $Z_2$ we have
\begin{theorem}\label{thmsm1}All the twisted spreading models of $Z(\mathcal T^2)$ are isomorphic to $\ell_2$.
\end{theorem}
To prove the result we need a variation on the inequality \cite[Lemma 4.8.]{CaFG} for almost optimal selections and we have to be delicate with the estimates. Therefore, let us recall that given a block $u$ with $\|u\|_2=1$, we may split it using the Lozanovskii factorization as 
\begin{equation}\label{split}
|u|^2=|v|\cdot|w|,
\end{equation}
where $\|v\|_X,\|w\|_{X^*}\leq 1$ and $\supp u=\supp v=\supp w$, see e.g. \cite{KM}. Assume $X, X^*$ have a joint 1-unconditional basis. Therefore, we have 
\begin{eqnarray*}
1=\sum_{j=1}^{\infty} |u(j)|^2&=&\sum_{j=1}^{\infty} |v(j)||w(j)|\\
&=& \langle |u|,|v| \rangle\\
&\leq& \|v\|\|w\|,
\end{eqnarray*}
 by the biorthogonality of $|v|=\sum_{j=1}^{\infty} |v(j)|e_j$ and $|w|=\sum_{j=1}^{\infty} |w(j)|e_j^*$. So we find $1\leq \|v\|^{1/2}\|w\|^{1/2}\leq 1$, and thus
\begin{equation*}
\|u\|_2=\|v\|_X=\|w\|_{X^*}=1.
\end{equation*}
We shall deal with selectors for the quotient map $\delta_{1/2}:\mathcal F_{\infty}(X,X^*)\to \ell_2$. In this sense, let us observe that the map
$$Bu(z):=u|v|^{1/2-z}|w|^{z-1/2}\in \mathcal F_{\infty}(X,X^*),$$
is a well defined selector for $u$ using the choices above. In particular, for every $t\in \mathbb R$
$$\|Bu(0+it)\|_X=\||u||v|^{1/2}|w|^{-1/2}\|_X=\||v|\|_X=1,$$ where both equalities hold by the $1$-unconditionality of the basis of $X$. The same arguments show that $\|Bu(1+it)\|_{X^*}=1$ for every $t\in \mathbb R$. Once these small technicalities have been clarified, let us state our spreading version of the inequality \cite[Lemma 4.8.]{CaFG}.
\begin{lemma}\label{estimate1}
Given $(u_j)_{j=1}^{\infty}$, let $|u_j|^2=|v_j||w_j|$ as in (\ref{split}) for $\mathcal T^2$ and $(\mathcal T^2)^*$, where $j\in \mathbb N$. If $\Omega_2$ denotes the centralizer induced by $(\mathcal T^2, (\mathcal T^2)^*)_{1/2}$, then there is $C>0$ such that for any sequence $b_1,...,b_n$ of scalars we have that
\begin{equation}\label{eq1}
\left\| \Omega_2\left( \sum_{j=1}^n b_j u_j\right)-\sum_{j=1}^n b_j \Omega_2(u_j)-\log \frac{\left\| \sum_{j=1}^n b_jv_j \right\|_{\mathcal T^2}}{\left\| \sum_{j=1}^n b_jw_j \right\|_{(\mathcal T^2)^*}}\sum_{j=1}^n b_ju_j\right\|_2
\end{equation} is $C$-bounded by
$$\left\|\sum_{j=1}^n b_ju_j\right \|_2+ \left\| \sum_{j=1}^n b_jv_j \right\|_{\mathcal T^2}^{1/2} \left\| \sum_{j=1}^n b_jw_j \right\|_{(\mathcal T^2)^*}^{1/2}. $$
\end{lemma}
\begin{proof}
Let $$F(z):=\frac{\sum_{j=1}^n b_jBu_j(z)}{\|\sum_{j=1}^n b_jv_j\|_{\mathcal T^2}^{1-z}\|\sum_{j=1}^n b_jw_j\|_{(\mathcal T^2)^*}^{z}}\in \mathcal F_{\infty}(\mathcal T^2, (\mathcal T^2)^*).$$
The critical step to check that $F$ is well defined is just to show that there is an independent numerical constant $C$ for which $\|F\|_{\mathcal F_{\infty}}\leq C$. By our previous discussion, $\|Bu_j(0+it)\|=1$ for every $t\in \mathbb R$ and since trivially $\supp Bu_j(0+it)=\supp v_j$, we find that 
$$\frac{\left\|\sum_{j=1}^n b_jBu_j(0+it)\right\|_{\mathcal T^2}}{\left\|\sum_{j=1}^n b_jv_j\right\|_{\mathcal T^2}}\sim C_1,$$
by a simple application of \cite[Proposition X.e.2(a)]{CS}. The same holds for the blocks $Bu_j(1+it)$ with perhaps a different constant $C_2$. We are basically done since the constant $C=\max\{C_1,C_2\}$ is essentially the constant of the lemma. Let us see how to conclude, so observe that
\begin{equation*}
\left\| \sum_{j=1}^n b_jv_j \right\|_{\mathcal T^2}^{1/2} \left\| \sum_{j=1}^n b_jw_j \right\|_{(\mathcal T^2)^*}^{1/2}F'(1/2)
\end{equation*}
is exactly equal to
\begin{equation*}
\sum_{j=1}^n b_j \Omega_2(u_j)-\log \frac{\left\| \sum_{j=1}^n b_jv_j \right\|_{\mathcal T^2}}{\left\| \sum_{j=1}^n b_jw_j \right\|_{(\mathcal T^2)^*}}\sum_{j=1}^n b_ju_j.
\end{equation*}
Therefore, the quantity of (\ref{eq1}), say $\Lambda$, becomes trivially 
$$\Lambda=\left\| \delta_{1/2}'\left( B(\sum_{j=1}^n b_ju_j)-  \left\| \sum_{j=1}^n b_jv_j \right\|_{\mathcal T^2}^{1/2} \left\| \sum_{j=1}^n b_jw_j \right\|_{(\mathcal T^2)^*}^{1/2}F(z) \right) \right\|.$$
The inner function above belongs now to $\Ker \delta_{1/2}$, this is, we have a cancellation. In particular, using \cite[Theorem 4.1.]{CCS}, there is an absolute numerical constant $C_3$ so that 
\begin{eqnarray*}
\Lambda&\leq& C_3\left( 2\left\|\sum_{j=1}^n b_ju_j\right\|_2+  \left\| \sum_{j=1}^n b_jv_j \right\|_{\mathcal T^2}^{1/2} \left\| \sum_{j=1}^n b_jw_j \right\|_{(\mathcal T^2)^*}^{1/2} \|F\|_{\mathcal F_{\infty}} \right)\\
&\leq& 2C_3\cdot C \left( \left\|\sum_{j=1}^n b_ju_j\right\|_2+  \left\| \sum_{j=1}^n b_jv_j \right\|_{\mathcal T^2}^{1/2} \left\| \sum_{j=1}^n b_jw_j \right\|_{(\mathcal T^2)^*}^{1/2}\right). 
\end{eqnarray*}
\end{proof}
We are ready to prove our claim on asymptotic models.
\begin{proof}[Proof of Theorem \ref{thmsm1}]
Let $\{(u_j,0),(\Omega_2(u_j),u_j)\}$ be, for simplicity, the sequence of the spreading. So we have to deal with an expression of the form $$\Lambda_{\nu}:=\left\| \sum_{j=1}^n a_ju_{k_j}+\sum_{j=1}^n b_j\Omega_2(u_{k_j})-\Omega_2\left(\sum_{j=1}^n b_ju_{k_j}\right) \right\|_2+\left\| \sum_{j=1}^n b_ju_{k_j}\right\|_2,$$
for $\nu\leq k_1<...< k_n$. Split each $u_{k_j}$ as in (\ref{split}) so that $|u_{k_j}|^2=|v_{k_j}||w_{k_j}|$ with $\|v_{k_j}\|_{\mathcal T^2}=1=\|w_{k_j}\|_{(\mathcal T^2)^*}$. If $n<\nu$ we have $n$-normalized vectors supported after $n$, and thus 
\begin{equation}\label{equivalence}
\left\| \sum_{j=1}^n b_jv_{k_j} \right\|_{\mathcal T^2}\sim \left(\sum_{j=1}^n |b_j|^2\right)^{1/2}\sim \left\| \sum_{j=1}^n b_jw_{k_j} \right\|_{(\mathcal T^2)^*}.
\end{equation}
Therefore, using Lemma \ref{estimate1}, we easily find, adding and extracting the log factor, that for $n<\nu\leq k_1<...<k_n$
\begin{eqnarray*}\label{bounded}
\Lambda_{\nu}&\leq&\left\| \sum_{j=1}^n a_ju_{k_j} \right\|_2+C_1 \left\|\sum_{j=1}^n  b_ju_{k_j} \right\|_2\\
&\leq& C_2\|(a,b)\|_{\ell_2\oplus_2\ell_2},
\end{eqnarray*}
where the constant $C_1$ depends only on the constant $C$ of Lemma \ref{estimate1} and the constant of the equivalence (\ref{equivalence}). Taking limit as $\nu \to \infty$, we get one inequality. The reverse inequality is proved in a similar vein but we include the argument for the sake of completeness. In a similar vein as before, we have
\begin{eqnarray*}\label{bounded2}
\left\| \sum_{j=1}^n a_ju_{k_j} \right\|_2+ \left\|\sum_{j=1}^n  b_ju_{k_j} \right\|_2&\leq&\Lambda_{\nu}+\left\| \sum_{j=1}^n b_j\Omega_2(u_{k_j})-\Omega_2\left(\sum_{j=1}^n b_ju_{k_j}\right) \right\|_2\\
&\leq& \Lambda_{\nu}+C\cdot\Delta(b,u)\\
&+&\left\|\log \frac{\left\| \sum_{j=1}^n b_jv_{k_j} \right\|_{\mathcal T^2}}{\left\| \sum_{j=1}^n b_jw_{k_j} \right\|_{(\mathcal T^2)^*}}\sum_{j=1}^n b_ju_{k_j}\right\|_2\\
&\leq&C_1 \cdot \Lambda_{\nu},
\end{eqnarray*}
where the first inequality follows if we add and extract the appropriate factors while the second employs Lemma \ref{estimate1}. The last inequality holds trivially if the blocks $u_{k_j}$ are chosen down far enough which is the case as $\nu \to \infty$.
\end{proof}
The same technique proves the following.
\begin{proposition}\label{smzt2block}
\label{spreadinormal}
All the spreading models of $Z(\mathcal T^2)$ are isomorphic to $\ell_2$.
\end{proposition}
\begin{proof} Let $(w_j)$ be a normalized and weakly null sequence in $Z(\mathcal T^2)$. There is no loss of generality assuming that it is a block basic sequence. Thus, let $(x_j,u_j)$ be such sequence. If $\|u_j\|\to 0$ we find easily a subsequence equivalent to the unit vector basis of $\ell_2$ since $$\left \|(x_j,u_j)-(x_j-\Omega_2(u_j),0)\right\|=\left\|(\Omega_2(u_j),u_j) \right\|=\|u_j\|.$$ So let us assume that there is $\eta >0$ and a further subsequence (we do not relabel) so that $\|u_j\|\geq \eta$. Now, it is easy to check that for every choice of scalars $b_1,...,b_n$, the expressions
$$\left\|\sum_{j=1}^nb_j\cdot x_j-\Omega_2\left (\sum_{j=1}^nb_j u_j\right)\right\|+\left\|\sum_{j=1}^nb_j u_j\right\|,$$ 
and 
$$\left\|\sum_{j=1}^nb_j\cdot \Omega_2(u_j)-\Omega_2\left (\sum_{j=1}^nb_j u_j\right)\right\|+\left\|\sum_{j=1}^nb_ju_j\right\|,$$ 
are equivalent since $\|x_j-\Omega_2(u_j)\|\leq 1$. Using the fact that the blocks $u_j$ are semi-normalized, we may use the arguments of Theorem \ref{thmsm1} with $a_1=...=a_n=0$ to conclude. 
\end{proof}
\section{The Principle of Small Perturbations for $Z(\mathcal T^2)$}\label{sectionfz2}
The claim on the P.S.P. for $Z(\mathcal T^2)$ is one of the main results of the paper. The key step to prove the P.S.P. for $Z(\mathcal T^2)$ is to read Lemma \ref{estimate1} in a different way, not a spreading but an optimal way.
\begin{theorem}\label{thmsingzt2}
The space $Z(\mathcal T^2)$ satisfies The Principle of Small Perturbations.
\end{theorem}
\begin{proof}
Let $(w_j)$ be a sequence in $Z(\mathcal T^2)$ which is equivalent to the unit vector basis of $\ell_2$ and, since $w_j\to 0$ weakly, pass to a subsequence which we do not relabel so that there exists semi-normalized blocks $(x_j,u_j)$ with $$\|w_j-(x_j,u_j)\|\leq 2^{-j},\;\;\;j\in \mathbb N.$$ We claim that $\|u_j\|\to 0$ and the rest of the proof will deal with this. So we assume that this does not hold and there is a further subsequence which again we do not relabel so that $\|u_j\|\geq \eta>0$. Since the blocks are semi-normalized, there is $C_1>0$ so that for every $j\in \mathbb N$
$$\|x_j-\Omega_2(u_j)\|+\|u_j\|\leq C_1.$$
To avoid entering more constants, we shall assume $\|u_j\|=1$ for all $j\in \mathbb N$. Thus, to be equivalent to the unit vector basis of $\ell_2$ means that for every $a_1,...,a_n$ we have
$$\left\| \sum_{j=1}^na_jx_j-\Omega_2\left(\sum_{j=1}^n a_j u_j \right) \right\|+ \left\| \sum_{j=1}^n a_j u_j \right\|\leq C_2\left \|\sum_{j=1}^n a_j u_j \right\|.$$
In particular, since $\|x_j-\Omega_2(u_j)\|\leq C_1$, we find that for every $a_1,...,a_n$
\begin{equation}\label{contr}
\left\| \sum_{j=1}^na_j\Omega_2(u_j)-\Omega_2\left(\sum_{j=1}^n a_j u_j \right) \right\|\leq C_3\left\|\sum_{j=1}^n a_j u_j \right\|.
\end{equation}
Let us check that (\ref{contr}) cannot hold. Pick $|u_j|^2=|v_j|\cdot |w_j|$ as in (\ref{split}) and write $$v_j=\sum_{k=p_j+1}^{p_{j+1}}\alpha_{k} t_{k}, \;\;\;\;w_j=\sum_{k=p_j+1}^{p_{j+1}}\beta_{k} t_{k}^*.$$
By \cite[Proposition X.e.2(a)]{CS}, we have that  $$\left \|\sum_{j=1}^{\infty}b_j v_j \right\|_{\mathcal T^2} \sim \left \|\sum_{j=1}^{\infty} b_j t_{p_{j}+1} \right\|_{\mathcal T^2}, $$
and also trivially that
$$\left \|\sum_{j=1}^{\infty}b_j w_j \right\|_{(\mathcal T^2)^*} \sim \left \|\sum_{j=1}^{\infty} b_j t_{p_{j}+1}^* \right\|_{(\mathcal T^2)^*}, $$
for every $(b_j)_{j=1}^{\infty}\in c_{00}$. Thus, we may use Lemma \ref{estimate1} and deduce that
\begin{equation*}\label{eq2}
\left\| \Omega_2\left( \sum_{j=1}^n b_j u_j\right)-\sum_{j=1}^n b_j \Omega_2(u_j)-\log \frac{\left \|\sum_{j=1}^{n} b_j t_{p_{j}+1} \right\|_{\mathcal T^2}}{\left \|\sum_{j=1}^{n} b_jt_{p_{j}+1}^* \right\|_{(\mathcal T^2)^*}}\sum_{j=1}^n b_ju_j\right\|_2
\end{equation*} 
is $C$-bounded by 
$$\left\|\sum_{j=1}^n b_ju_j\right \|_2+ C_4\left \|\sum_{j=1}^{n} b_j t_{p_{j}+1} \right\|_{\mathcal T^2}^{1/2}\left \|\sum_{j=1}^{n} b_j t_{p_{j}+1}^* \right\|_{(\mathcal T^2)^*}^{1/2},$$ for perhaps a different $C>0$ and where $C_4$ is independent of $n\in \mathbb N$. Let us pick a very special choice of $(b_j)_{j=1}^n$. To do it, we use the method of critical points developed in \cite{JS4}. Thus, we first define $$\kappa(n):=\sup \left\{ \left\|\sum_{j=1}^n a_je_{p_{j}+1} \right\|_{\ell_2}: \left\|\sum_{j=1}^n a_jt_{p_{j}+1} \right\|_{\mathcal T^2}=1 \right\},$$
which is attained by compactness. Pick $0\neq a=(a_1,...,a_n)$ so that 
\begin{equation}\label{extremal}\left\|\sum_{j=1}^n a_je_{p_{j}+1} \right\|_{\ell_2}=\kappa(n)\left\|\sum_{j=1}^n a_jt_{p_{j}+1} \right\|_{\mathcal T^2}.
\end{equation}
Since the basis of $\ell_2$ dominates the basis of $\mathcal T^2$ and thus all its subsequences, we find by simple duality that for every $b$ with $\supp b \subseteq \{p_{1}+1,...,p_n+1 \}$, we have 
\begin{equation}\label{extremal2}
\|b\|_{\ell_2}\leq \|b\|_{(\mathcal T^2)^*}\leq \kappa(n)\|b\|_{\ell_2}.
\end{equation}
In particular, the inequalities above hold for $a=\sum_{j=1}^n a_jt_{p_{j}+1}^*$. Observe that for our choice $a$, our version of the lemma gives that 
\begin{eqnarray*}
\Delta(a,u)&=&\left\|\sum_{j=1}^n a_ju_j\right \|_2+ C_4\left \|\sum_{j=1}^{n} a_j t_{p_{j}+1} \right\|_{\mathcal T^2}^{1/2}\left \|\sum_{j=1}^{n} a_j t_{p_{j}+1}^* \right\|_{(\mathcal T^2)^*}^{1/2}\\
&=& \left\|\sum_{j=1}^n a_ju_j\right \|_2+ C_4\kappa(n)^{-1/2}\left \|\sum_{j=1}^{n} a_j t_{p_{j}+1} \right\|_{2}^{1/2}\left \|\sum_{j=1}^{n} a_j t_{p_{j}+1}^* \right\|_{(\mathcal T^2)^*}^{1/2}\\
&\leq&(C_4+1)\left\|\sum_{j=1}^n a_ju_j\right \|_2.
\end{eqnarray*}
We are about to see that we have a contradiction. If we plug (\ref{extremal}) in (\ref{extremal2}) we have trivially $$\kappa(n)\left \|\sum_{j=1}^{n} a_j t_{p_{j}+1} \right\|_{\mathcal T^2}\leq \left \|\sum_{j=1}^{n} a_j t_{p_{j}+1}^* \right\|_{(\mathcal T^2)^*} \leq \kappa(n)^2 \left \|\sum_{j=1}^{n} a_j t_{p_{j}+1} \right\|_{\mathcal T^2},$$ and in particular
$$\log \kappa(n)\sim \left|\log \frac{\left \|\sum_{j=1}^{n} a_j t_{p_{j}+1} \right\|_{\mathcal T^2}}{\left \|\sum_{j=1}^{n} a_j t_{p_{j}+1}^* \right\|_{(\mathcal T^2)^*}}\right|.$$
We use the previous estimate and the special form of the lemma described above to finally have
\begin{eqnarray*}
\log \kappa(n) \left\| \sum_{j=1}^n a_ju_j\right\|_2 &\leq&C\cdot \Delta(a,u) +\left\| \sum_{j=1}^na_j\Omega_2(u_j)-\Omega_2\left(\sum_{j=1}^n a_j u_j \right) \right\|\\
&\leq&(C(C_4+1)+C_3)\left\|\sum_{j=1}^n a_ju_j\right \|_2,
\end{eqnarray*}
where we have plugged (\ref{contr}) to bound with $C_3$ in the second inequality.
Therefore, $\log \kappa(n)$ and thus $\kappa(n)$ is uniformly bounded for every $n\in \mathbb N$. This proves that $\{t_{p_{j}+1}\}_{j=1}^{\infty}$ spans an isomorphic copy of $\ell_2$ in $\mathcal T^2$ which is absurd. So it must be $\|u_j\|\to 0$ and the claim is proved. We may argue safely now as in the proof of Proposition \ref{PSP}: Let us assume $\|u_j\|_2\leq 2^{-j}$, and hence $$\|w_j-(y_j,0)\|_2\leq 2^{-j+1},$$
where $y_j=x_j-\Omega(u_j)$. Replacing $w_j$ by $\lambda_j w_j$, where the sequence $\lambda_j$ satisfies $\inf \lambda_j>0$ and $\sup_j \lambda_j<\infty$, we may assume that $\|y_j\|_2=1$.
\end{proof}
We find easily now a result on singularity.
\begin{corollary}\label{zt2ss}
 The natural quotient of $Z(\mathcal T^2)$ onto $\ell_2$ is strictly singular.
\end{corollary}
\begin{proof}
The claim is similar to Corollary \ref{one} which only uses the P.S.P.
\end{proof}
Since it is well known that $Z(\ell_2)=\ell_2$, one trivially has an example of an asymptotic $\ell_2$-space whose derivation fails to have a strictly singular quotient map onto $\ell_2$. Therefore, the previous corollary shows an extreme situation and a dichotomy for $p=2$. It also completes the cases $p\neq 2$ of Castillo \textit{et al}. in \cite{CaFG} that we already obtained in Corollary \ref{one}. Observe that in our approach all cases are unified as a consequence of the P.S.P.
\section{The spreading models of $Z(\mathcal T)$}
We compute the asymptotic models of $Z(\mathcal T)$ using a similar technique as in the previous section. For simplicity we use a compact notation, so we do not write sums or subscripts. For example, for $u=(u_1,...,u_n)$ and $b=(b_1,...,b_n)$ we write $$b\cdot u=\sum_{j=1}^n b_ju_j.$$
With this notation, let us state a technical lemma as in the previous section. 
\begin{lemma}\label{estimate2}
Given $(u_j)_{j=1}^{\infty}$, let $|u_j|^2=|v_j||w_j|$ as in (\ref{split}) for $\mathcal T$ and $\mathcal T^*$, where $j\in \mathbb N$. If $\Omega_1$ denotes the centralizer induced by $(\mathcal T, \mathcal T^*)_{1/2}$, then there is $C>0$ such that for any sequence $b_1,...,b_n$ of scalars we have
\begin{equation*}\label{eq3}
\left\| \Omega_1\left( b\cdot u\right)-b\cdot \Omega_1(u)-\log \frac{\left\| |b|^2\cdot v \right\|_{\mathcal T}}{\left\|  b|b|^{-1}\cdot w \right\|_{\mathcal T^*}}b\cdot u+2b\log|b|\cdot u\right\|_2\leq C\cdot \Delta(b,u),
\end{equation*} where
$$\Delta(b,u)=\left\|b\cdot u\right \|_2+ \left\| |b|^2\cdot v \right\|_{\mathcal T}^{1/2} \left\| b|b|^{-1} \cdot w \right\|_{\mathcal T^*}^{1/2}. $$
\end{lemma}
\begin{proof}The proof is similar to that of Lemma \ref{estimate1} using the function $$G(z):=\frac{\sum_{j=1}^n b_j|b_j|^{-2(z-1/2)}Bu_j(z)}{\left \|\sum_{j=1}^n|b_j|^2v_j\right\|_{\mathcal T}^{1-z}\left \|\sum_{j=1}^nb_j|b_j|^{-1}w_j\right\|_{\mathcal T^*}^{z}}$$
with the aid of \cite[Proposition II.4]{CS} to have control on the norm of $G$.
\end{proof}
As a consequence, we may study the asymptotic structure of $Z(\mathcal T)$.
\begin{theorem}\label{thmsm2}
All the twisted spreading models of $Z(\mathcal T)$ are isomorphic to $Z_2$. 
\end{theorem}
\begin{proof}
Let $\mathcal K$ denotes the Kalton-Peck map arising from the couple $(\ell_1,c_0)$ and let $u=(u_{k_1},...,u_{k_n})$. Using the triangle inequality for $$\nabla (a,b):=\left\|a\cdot u- \Omega_1(b\cdot u) +b \cdot\Omega_1(u)\right\|,$$ we have
\begin{eqnarray*}
\nabla (a,b)&\leq& \left\|a\cdot u-\mathcal K(b\cdot u) +b \cdot\mathcal K(u)\right\|\\
&+& \left\|\mathcal K(b\cdot u) -b \cdot\mathcal K(u)-\log\frac{\||b|^2\cdot v\|_{\mathcal T}}{\|b|b|^{-1} \cdot w\|_{\mathcal T^*}}b\cdot u  +2b\log|b|u\right\|\\
&+&\left\|-\Omega_1(b\cdot u) +b \cdot\Omega_1(u)+\log\frac{\||b|^2\cdot v\|_{\mathcal T}}{\|b|b|^{-1} \cdot w\|_{\mathcal{T}^*}}b\cdot u -2b\log|b|u\right\|.
\end{eqnarray*}
For $n<\nu\leq k_1<...<k_n$ we have $n$ normalized vectors supported after $n$, so that $$\||b|^2\cdot v\|_{\mathcal T}\sim \sum_{j=1}^n |b_j|^2,$$ and $$\|b|b|^{-1} \cdot w\|_{\mathcal T^*}\sim 1.$$ This is, at the level $n<\nu$ we may not distinguish, up to a numerical constant, the $\ell_1$-norm from the $\mathcal T$-norm of $n$ normalized blocks (and similarly with $c_0$ and $\mathcal T^*$). Let us observe that Lemma \ref{estimate2} also holds for $\mathcal K$ replacing $\mathcal T$ and $\mathcal T^*$ by $\ell_1$ and $c_0$ respectively and hence  
$$ \left\|\mathcal K(b\cdot u) -b \cdot\mathcal K(u)-\log\frac{\||b|^2\cdot v\|_{\ell_1}}{\|b|b|^{-1} \cdot w\|_{c_0}}b\cdot u  +2b\log|b|u\right\|\leq C \|b\cdot u \|_2.$$
 Then, by the previous comments, we find that for some $C_1>0$
 \begin{equation}\label{zk1}
 \left\|\mathcal K(b\cdot u) -b \cdot\mathcal K(u)-\log\left(\frac{\||b|^2\cdot v\|_\mathcal{T}}{\|b|b|^{-1} \cdot w\|_{\mathcal{T}^*}\cdot|b|^{2}}\right) b\cdot u  \right\|\leq C_1 \|b\cdot u \|_2.
\end{equation}
If we apply Lemma \ref{estimate2} to $\Omega_1$, we find a $C_2>0$ such that
 \begin{equation}\label{zk2}
\left\|\Omega_1(b\cdot u) -b \cdot\Omega_1(u)-\log\left(\frac{\||b|^2\cdot v\|_\mathcal{T}}{\|b|b|^{-1} \cdot w\|_{\mathcal{T}^*}\cdot|b|^{2}}\right)b\cdot u \right\|\leq C_2  \|b\cdot u \|_2.
\end{equation}
If we apply the estimates (\ref{zk1}) and (\ref{zk2}) to the first inequality on $\nabla (a,b)$, we will find a $C_3$ so that 
\begin{equation}\label{zk3}
\nabla (a,b)+\|b\cdot u\|\leq C_3\left(\left\|a\cdot u-\mathcal K(b\cdot u) +b \cdot\mathcal K(u)\right\|+\|b\cdot u\|\right),
\end{equation}
whenever $n<\nu\leq k_1<...<k_n$ for $u=(u_{k_1},...,u_{k_n})$ which is the case as $\nu \to \infty$. The reverse inequality is completely similar choosing blocks starting down far enough.
\end{proof}
 We have now also the analogous result of Example \ref{example2} for the spreading models of $Z(\mathcal T)$.
\begin{proposition}\label{2spr} The spreading models of $Z(\mathcal T)$ are the Hilbert copy $\ell_2$ and the Orlicz space $\ell_M$. 
\end{proposition}
\begin{proof} There is no loss of generality assuming the sequence is given by blocks. So, let $(x_j,u_j)$ be such sequence of blocks. If $\|u_j\|\to 0$ we find easily a subsequence equivalent to the unit vector basis of $\ell_2$ arguing exactly as in Proposition \ref{spreadinormal}. So let us assume there is $\eta >0$ and a further subsequence which we do not relabel so that $\|u_j\|\geq \eta$. Now, it is easy to check that for every choice of scalars $b_1,...,b_n$, the expressions
$$\|b\cdot x-\Omega_1(b\cdot u)\|+\|b\cdot u\|,$$ 
and 
$$\|b\cdot \Omega_1(u)-\Omega_1(b\cdot u)\|+\|b\cdot u\|$$
are equivalent since $\|x_j-\Omega(u_j)\|\leq 1$. Using the fact that the blocks $u_j$ are semi-normalized, we may use inequality (\ref{zk3}) with $a=0$ and its reversal to show that both expressions are indeed equivalent to 
$$\|b\cdot \mathcal K(u)-\mathcal K(b\cdot u)\|+\|b\cdot u\|,$$
whenever the blocks are chosen down far enough. As we commented troughout Example \ref{example}, with the aid of \cite[Theorem 6.5, page 28, line 17]{KaPe}, the expression above is equivalent to the norm of $\|(0,b)\|_{Z_2}$, so we are done. 
\end{proof}

We checked in Corollary \ref{maxmin} that $Z(\mathcal T)$ does not embed into $Z(\mathcal T^2)$. Now, as an application of the previous proposition, we distinguish $Z(\mathcal T^2)$ also from a subspace of $Z(\mathcal T)$.
\begin{proposition}\label{zt2zt}
The twisted Hilbert space $Z(\mathcal T^2)$ is not isomorphic to a subspace or a quotient of $Z(\mathcal T)$.
\end{proposition}
\begin{proof}
Assume that there is such an embedding $T:Z(\mathcal T^2)\to Z(\mathcal T)$. Passing to a subsequence if necessary we may assume that the sequence $\{T(0,e_j)\}\jinf$ is equivalent to a block basic sequence $\{(x_j,u_j)\}\jinf$. If $\|u_j\|\to 0$ there is a further infinite subset $\mathbb M\subseteq \mathbb N$ so that $\{(0,e_j)\}_{j\in \mathbb M}$ is equivalent to unit vector basis of $\ell_2$. This cannot hold by a simple application of Kalton uniqueness theorem \cite[Theorem 7.6.]{ka}. Indeed, we would find that for a suitable $C>0$ and every $(b_j)_{j\in \mathbb M}$ $$\left\| \sum_{j\in \mathbb M} b_j (0,e_{j})\right\|=\left \|\Omega_2\left(\sum_{j\in \mathbb M}b_j e_j\right)\right \|+ \left \| \sum_{j\in \mathbb M}b_j e_j\right \|\leq C \left\| \sum_{j\in \mathbb M}b_j e_j\right \|.$$ In other words, the centralizer $\Omega_{2|\mathbb M}$ is $(C-1)$-bounded. Let us recall the interpolation formula $$(\mathcal T^2(\mathbb M), (\mathcal T^2)^*(\mathbb M))_{1/2}=\ell_2(\mathbb M)$$ (see \cite[Corollary 4.3]{CoS}) for which we may trivially assume that $\Omega_{|\mathbb M}$ is the corresponding centralizer. Since $\Omega_{2|\mathbb M}$ is bounded, we find that $\mathcal T^2(\mathbb M)\approx\ell_2(\mathbb M)$  by the aforementioned Kalton's result \cite[Theorem 7.6.]{ka} which is impossible since $\mathcal T^2$ contains no copies of $\ell_2$ \cite{CS}. So it must be that there exists a further subsequence so that $\|u_{k_j}\|\geq \eta$. But this sequence $\{(x_{k_j},u_{k_j})\}\jinf$ generates by Proposition \ref{2spr} a spreading model that is isomorphic to $\ell_M$. On the other hand, by Proposition \ref{smzt2block}, the only spreading model that generates $\{(0,e_{k_j})\}\jinf$ in $Z(\mathcal T^2)$ is $\ell_2$. So we get a contradiction. The claim for the quotient follows by duality since both spaces are isomorphic to its own dual.
\end{proof}
%
\section{Twisted Hilbert spaces induced by asymptotic $\ell_p$-spaces}
We show that there is a large number of twisted Hilbert spaces having $Z_2$ as a twisted spreading model. This is only the abstract statement of the ideas given in the previous sections. We begin with a technical lemma.
\begin{lemma}\label{estimate3}
Let $X$ be a space with a shrinking unconditional basis. Given $(u_j)_{j=1}^{\infty}$, let $|u_j|^2=|v_j||w_j|$ as in (\ref{split}) for $X$ and $X^*$, where $j\in \mathbb N$. If $\Omega_X$ denotes the centralizer induced by $(X,X^*)_{1/2}$ and we let $1\leq p,q\leq\infty$ with $p^{-1}+q^{-1}=1$, then there is $C>0$ such that for any sequence $b_1,...,b_n$ of scalars we have that
\begin{equation*}
\left\| \Omega_X\left( b\cdot u\right)-b\cdot \Omega_X(u)-\log \frac{\left\| |b|^{2/p}\cdot v \right\|_{X}}{\left\|  |b|^{2/q}\cdot w \right\|_{X^*}}b\cdot u-\left(\frac{2}{q}-\frac{2}{p}\right)b\log|b|\cdot u\right\|_2
\end{equation*} is $C$-bounded by 
$$\left\|b\cdot u\right \|_2+ \left\| |b|^{2/p}\cdot v \right\|_{X}^{1/2} \left\| |b|^{2/q} \cdot w \right\|_{X^*}^{1/2}. $$
\end{lemma}
\begin{proof}The proof is similar to that of Lemma \ref{estimate1} using the function $$G(z):=\frac{\sum_{j=1}^n b_j|b_j|^{\left(\frac{2}{q}-\frac{2}{p}\right)(z-1/2)}Bu_j(z)}{\left \|\sum_{j=1}^n|b_j|^{2/p}v_j\right\|_{X}^{1-z}\left \|\sum_{j=1}^n|b_j|^{2/q}w_j\right\|_{X^*}^{z}}.$$
The bound for the norm of $G$ follows by the unconditionality of the basis of $X,X^*$ and the special form of the map $B$ described after (\ref{split}).
\end{proof}
Now, we give the result.
\begin{theorem}\label{thmrefl}
Let $X$ be a space with a shrinking unconditional basis. If $X$ is an asymptotic $\ell_p$-space with $1\leq p\neq2<\infty$, then all the twisted spreading models of $Z(X)$ are isomorphic to the Kalton-Peck space $Z_2$.
\end{theorem}
\begin{proof}
The proof is similar to that of Theorem \ref{thmsm2} using the lemma above. We only need to recall  a couple of facts. First, under the assumptions of the theorem, $X^*$ is an asymptotic $\ell_q$-space with $1/p+1/q=1$. And second, the Kalton-Peck map $\mathcal K$ of Theorem \ref{thmsm2} is the one arising from the couple $(\ell_1, c_0)$. For our case we need the Kalton-Peck map arising from the couple $(\ell_p,\ell_q)$.
\end{proof}
For the case of $c_0$, the following proposition shows how to find examples of such phenomena using only spreading models.
\begin{proposition}\label{propasym}
Let $X$ be a space with a shrinking unconditional basis $(e_j)_{j=1}^{\infty}$. Assume that every spreading model of a block basic sequence of $(e_j)_{j=1}^{\infty}$ is isomorphic to $c_0$. Then, all the twisted spreading models of $Z(X)$ are isomorphic to the Kalton-Peck space $Z_2$.
\end{proposition}
\begin{proof}
Let $\{(u_j,0),(\Omega(u_j),u_j)\}_{j=1}^{\infty}$ be an spreading sequence in $Z(X)$, where $(u_j)_{j=1}^{\infty}$ is a normalized block basic sequence in $\ell_2$. Split $|u_j|^2=|v_j|\cdot |w_j|$ as in (\ref{split}), where $\|v_j\|_X=1, \|w_j\|_{X^*}=1$ for each $j\in \mathbb N$. Passing to further subsequences, we may assume that $(v_j)$ generates some spreading model. By assumption, there is some $C$ so that the spreading model of $(v_j)$ is $C$-isomorphic to $c_0$.  By the shrinking condition $(v_j)$ is weakly null, and thus the basis of the spreading is unconditional. Therefore, such basis is equivalent to the canonical basis of $c_0$. The key point is that the spreading model of $(|w_j|)$ is $\ell_1$. Indeed, for any $(\theta_j)_{j=1}^{\infty}$ with $|\theta_j|=1$ for each $j\in \mathbb N$, we have
\begin{eqnarray*}
\sum_{j=1}^n b_j^2=\left\|\sum_{j=1}^n b_ju_{k_j}\right\|_2^2&=&\sum_{j=1}^nb_j^2|u_{k_j}|^2\\
&=&\sum_{j=1}^nb_j^2|v_{k_j}||u_{k_j}|\\
&=&\sum_{j=1}^n \theta_jb_j^2\theta_j|v_{k_j}||w_{k_j}|=(*).\\
\end{eqnarray*}
Using the argument after (\ref{split}), we find the upper bound
\begin{eqnarray*}
(*)&\leq&\left\|\sum_{j=1}^n\theta_j|v_{k_j}|  \right\|_X \left\|\sum_{j=1}^n \theta_jb_j^2|w_j| \right\|_{X^*}\\
&\leq& 2KC\left\|\sum_{j=1}^n \theta_jb_j^2|w_j| \right\|_{X^*},\\
\end{eqnarray*}
where the last inequality holds if $k_1<...<k_n$ is large enough since $(e_j)$ is $K$-unconditional and thus the spreading of $(|v_j|)$ is also $KC$-isomorphic to $c_0$. Indeed, $$\left\|\sum_{j=1}^n\theta_j|v_{k_j}|  \right\|_X \leq K\left\|\sum_{j=1}^n\theta_jv_{k_j}  \right\|_X.$$
Passing to a further subsequence, we may assume that the spreading model of $(|w_j|)$ is equivalent to the canonical basis of $\ell_1$. Therefore, it is plain to compute that for a fixed $b_1,...,b_n$, letting $k_1<...<k_n$ large enough, the map $$G(z)=\frac{\sum_{j=1}^n b_j|b_j|^{2(z-1/2)}Bu_j(z)} {\left\|\sum_{j=1}^nb_j|b_j|^{-1}|v_{k_j}| \right\|_X^{1-z}\left\|\sum_{j=1}^n |b_j|^2|w_{k_j}| \right\|_{X^*}^z}\in \mathcal F_{\infty}(X,X^*)$$
is well defined and satisfies $$\|G\|_{\mathcal F_{\infty}}\leq C_1,$$
for some absolute $C_1$. If we denote $\Omega_X$ the corresponding centralizer, we find using the cancellation principle of Theorem \ref{thmsm2} that the quantity
\begin{equation*}\resizebox{\hsize}{!}{$
\left\|\Omega_X(\sum_{j=1}^nb_ju_{k_j})-\sum_{j=1}^nb_j\Omega_X(u_{k_j})-\sum_{j=1}^n b_j\log \left(\frac{\|\sum_{j=1}^nb_j|b_j|^{-1}|v_{k_j}|\|_X\cdot |b_j|^2}{\|\sum_{j=1}^n|b_j|^2w_{k_j}|\|_{X^*}}\right)u_{k_j}\right\|$}
\end{equation*}
is upper bounded, up to a constant, by 
$$\left\|\sum_{j=1}^nb_ju_{k_j}\right\|+\left\|\sum_{j=1}^nb_j|b_j|^{-1}|v_{k_j}|\right\|_X^{1/2}\left\|\sum_{j=1}^nb_j|b_j||w_{k_j}|\right\|_{X^*}^{1/2}.$$
Adding and extracting the factor $$\log \frac{\|\sum_{j=1}^nb_j|b_j|^{-1}e_j\|_{c_0}}{\|\sum_{j=1}^n b_j|b_j|e_j\|_{\ell_1}}\sum_{j=1}^nb_ju_{k_j},$$
in the expression of $\Omega_X$ above, we may argue exactly as in Theorem \ref{thmsm2} letting $k_1\to \infty$.
\end{proof}
For the next two corollaries, the space $X$ is in the hypotesis of Theorem \ref{thmrefl} or Proposition \ref{propasym}.
\begin{corollary}\label{smabstract}
The only spreading models of $Z(X)$ are the Hilbert copy $\ell_2$ and the Orlicz space $\ell_M$.
\end{corollary}
Thus, also
\begin{corollary}\label{smabstract2}
The twisted Hilbert space $Z(\mathcal T^2)$ is not isomorphic to a subspace or a quotient of $Z(X)$.
\end{corollary}

%


\section{The Kalton-Peck space as a spreading model}
We have seen roughly that if $X$ behaves asymptotically as $\ell_p$ or $c_0$ for $p\neq 2$, then $Z(X)$ looks asymptotically like $Z_2$. Let us go now in the opposite direction. Assume that a twisted Hilbert space $Z$ looks asymptotically as $Z_2$. What information can be recovered from $Z$? The key of our answer relies on the notion of maximal centralizer which closes our loop of ideas opened at Section \ref{sectionmM}.
\begin{theorem}\label{round}
Let $Z$ be a twisted Hilbert space with centralizer $\Omega$. If every twisted spreading model of $Z$ is isomorphic to the Kalton-Peck space $Z_2$, then $\Omega$ is maximal. In particular, $Z$ satisfies The Principle of Small Perturbations.
\end{theorem}
\begin{proof} Let $\Omega$ be the centralizer of $Z$ so let us prove first:\\

{\bf Claim M:} $\Omega$ is maximal.\\ 

Pick $(u_j)$ any sequence of semi-normalized blocks, say $\alpha\leq \|u_j\|\leq \beta$ for $j=1,2,...$, let $Z_S$ be the corresponding twisted spreading model (we do not relabel the subsequence) and $T:Z_S\longrightarrow Z_2$ the isomorphism. Thus, in particular, $$C^{-1}|(a,b)|_{Z_S}\leq \|T(a,b)\|_{Z_2}\leq C|(a,b)|_{Z_S},$$ for some fixed $C$ and every $(a,b)\in Z_S$. Let $x_j:=T(0,e_j)$ for $j\in \mathbb N$. There are two alternatives: either there is a subsequence such that $x_j$ is equivalent to the unit vector basis of $\ell_2$ or there is a subsequence such that the equivalence holds for the basis in the Orlicz space $\ell_M$, see \cite[Lemma 5.3.]{KaPe}. Since $T$ is an isomorphism, the same holds for a subsequence of $\{(0,e_j)\}_{j=1}^{\infty}$ and, by the invariance under spreading, it holds for the whole sequence $\{(0,e_j)\}_{j=1}^{\infty}$. Let us show the first possibility cannot hold. For if it is true then 
\begin{eqnarray*}
\left| \left(0,\sum_{j=1}^n b_je_j\right)\right|_{Z_S}&:=&\lims\left\|\Omega\left(\sum_{j=1}^n b_j u_{k_j} \right)-\sum_{j=1}^n b_j\Omega(u_{k_j}) \right\|\\
&+&\left\|\sum_{j=1}^nb_ju_{k_j}\right\|\\
\end{eqnarray*}
would be bounded by $$ C_1 \left(\sum_{j=1}^n |b_j|^2\right)^{1/2},$$
for some $C_1$. But arguing as in the proof of Theorem \ref{thmsm1}, we would find that the spreading model is isomorphic to $\ell_2$ which cannot be. Thus, it must equivalent to the Orlicz basis. In particular, there is $C_2$ such that for every $n\in \mathbb N$ we have
\begin{eqnarray*}
\lims\left\| \Omega\left(\sum_{j=1}^n u_{k_j} \right)-\sum_{j=1}^n \Omega(u_{k_j}) \right\| +\left\|\sum_{j=1}^nu_{k_j}\right\| &=&\left| \left( 0, \sum_{j=1}^ne_j\right)\right|_{Z_S}\\
&\geq& C_2\sqrt{n}\log{\sqrt{n}}.
\end{eqnarray*} 
Thus, by definition of limit, given $\varepsilon=1$ and $n\in \mathbb N$, there is $\nu\in \mathbb N$ such that 
$$\left\| \Omega\left(\sum_{j=1}^n u_{k_j} \right)-\sum_{j=1}^n \Omega(u_{k_j}) \right\| \geq \frac{C_2}{2}\sqrt{n}\log{\sqrt{n}}-\beta\sqrt{n},$$
whenever $\nu\leq k_1<...<k_n$. In other words, the subsequence of the spreading gives the good subsequence we are looking for. Since this can be done for every block basic sequence $(u_j)$, we find the maximality. So the claim is proved.\\

{\bf Claim P:} Maximal $\Rightarrow$ P.S.P.\\ 

This is done in its simplest form during the proof of Proposition \ref{PSP}. Let us squeeze the argument. Let $(w_j)$ be a sequence in $Z$ equivalent to the unit vector basis of $\ell_2$ and, since $w_j\to 0$ weakly, pass to a subsequence if necessary so that there exists blocks $(x_j,u_j)$ with $$\|w_j-(x_j,u_j)\|\leq 2^{-j},\;\;j\in \mathbb N.$$ We claim that $\|u_j\|\to 0$, so assume otherwise that there is further subsequence which we are not relabelling such that $\|u_j\|\geq \eta >0$. Then the blocks $(u_j)$ are semi-normalized so by maximality we may find a good subsequence, once again we do not relabel, say $(u_j)$. Let us observe that we have, for elements of this subsequence, the following crucial inequality
\begin{eqnarray*}\label{nidea}
\left\| \left(\sum_{j=1}^n x_{j},\sum_{j=1}^n u_{j}\right) \right\| &\geq& \left\| \sum_{j=1}^n x_{j}-\Omega \left( \sum_{j=1}^n u_{j} \right) \right\|\\
&\geq &  \left\| \Omega\left(\sum_{j=1}^n u_{j}\right)-\sum_{j=1}^n\Omega \left( u_{j} \right) \right\|-  \left\| \sum_{j=1}^n (x_{j}-\Omega \left(u_{j}\right) \right\|,
\end{eqnarray*}
which is impossible if $\|u_j\|\geq \eta >0$. Indeed, $\left\| \sum_{j=1}^n (x_j-\Omega \left(u_j\right) \right\|^2\leq \sum_{j=1}^n \|x_j-\Omega(u_j) \|^2\leq n$, while $\left\| \left(\sum_{j=1}^n x_j,\sum_{j=1}^n u_j\right) \right\| \sim \sqrt{n}$ because the blocks $(x_j,u_j)$ are a small perturbation of $(w_j)$. Now, by maximality, for each $n$ we may find $\nu$ such that for every $\nu\leq k_1<...<k_n$, we have 
$$\frac{\left\| \Omega\left(\sum_{j=1}^n u_{k_j}\right)-\sum_{j=1}^n\Omega \left( u_{k_j} \right) \right\|}{\sqrt{n}}\longrightarrow \infty.$$
Therefore, if for each $n\in \mathbb N$ we pick $k_1<...<k_n$ as above, it is impossible that the crucial inequality  holds for these choices if $n$ is large enough. Thus, $\|u_{j}\|\to 0$ if $j\to \infty$. The rest of the proof goes exactly as in Proposition \ref{PSP}. 
\end{proof}

\appendix


\begin{thebibliography}{99}

\bibitem{AK} F. Albiac and N.J. Kalton, \emph{Topics in Banach space theory}.
Graduate Texts in Mathematics 233. Springer-Verlag.

\bibitem{ACK} G. Androulakis, P.G. Casazza and D.N. Kutzarova, {\it Some more weak Hilbert spaces}, Canad. Math. Bull. 43 (2000), no. 3, 257--267.

%


\bibitem{BeLa} B. Beauzamy and J.T. Laprest\'e, \emph{Mod\`eles \'etal\'es des espaces de Banach}. (French) [Spreading models of Banach spaces] Publ. D\'ep. Math. (Lyon) (N.S.) 1983. 


\bibitem{BenLin} Y. Benyamini and J. Lindenstrauss, \emph{Geometric nonlinear functional analysis. Vol. 1}. American Mathematical Society Colloquium Publications, 48. American Mathematical Society, Providence, RI, 2000.



\bibitem{BeLo}  J. Bergh and J. L\"ofstr\"om, \emph{Interpolation spaces. An introduction}. Grundlehren der Mathematischen Wissenschaften, No. 223. Springer-Verlag, Berlin-New York, 1976.

\bibitem{BrSu} A. Brunel and L. Sucheston, {\it On B-convex Banach spaces}, Math. Systems Theory 7 (1974), no. 4, 294--299.




\bibitem{CS} P. Casazza \and T.J. Shura, \emph{Tsirelson's space. With an appendix by J. Baker, O. Slotterbeck and R. Aron}. Lecture Notes in Mathematics, 1363. Springer-Verlag, Berlin, 1989.

\bibitem{CN2} P.G. Casazza \and N.J. Nielsen, {\it A Banach space with a symmetric basis which is of weak cotype 2 but not of cotype 2}, Studia Math. 157 (2003), no. 1, 1--16.

\bibitem{CCS}  M.J. Carro, J. Cerd\`a  \and J. Soria, {\it Commutators and interpolation methods}, Ark. Mat. 33 (1995), no. 2, 199--216.


\bibitem{Ca} J.M.F. Castillo, {\it Personal communication}.

 

%


\bibitem{CaFG} J.M.F. Castillo, V. Ferenczi and M. Gonz\'alez, {\it Singular exact sequences generated by complex interpolation}, Trans. Amer. Math. Soc. 369 (2017), 4671-4708.

\bibitem{derivation} J.M.F. Castillo, D. Morales and J. Su\'arez de la Fuente, {\it Derivation of vector-valued complex interpolation scales}, J. Math. Anal. Appl. 468 (2018), no. 1, 461--472.

\bibitem{CoS} F. Cobos \and T. Schonbek, {\it On a theorem by Lions and Peetre about interpolation between a Banach space and its dual}, Houston J. Math. 24 (1998), no. 2, 325--344.


\bibitem{ELP} P. Enflo, J. Lindenstrauss \and G. Pisier, {\it On the ``three space problem"},  Math. Scand. 36 (1975), no. 2, 199--210.

\bibitem{FR} V. Ferenczi \and C. Rosendal, {\it Ergodic Banach spaces}, Adv. Math. 195 (2005), 259--282.

\bibitem{G} T. Gowers, {\it Gowers's Weblog}, https://gowers.wordpress.com/

\bibitem{GoMau} T. Gowers and B. Maurey, {\it The unconditional basic sequence problem}, J. Amer. Math. Soc. 6 (1993), no. 4, 851–-874. 

\bibitem{HalOdell} L. Halbeisen and E. Odell, {\it On asymptotic models in Banach spaces}, Israel J. Math. 139 (2004), 253--291.

\bibitem{Hal} P.R. Halmos, {\it I want to be a mathematician. An automatography.} Springer-Verlag, New York, 1985.

\bibitem{James} R.C. James, {\it Uniformly non-square Banach spaces}, Ann. of Math. (2) 80 (1964), 542–-550.



 \bibitem{ka} N.J. Kalton, {\it Differentials of complex interpolation processes for K\"othe function spaces}, Trans. Amer. Math. Soc. 333 (1992), no. 2, 479--529.
 
 \bibitem{ka2} N.J. Kalton, {\it Twisted Hilbert spaces and unconditional structure}, J. Inst. Math. Jussieu 2 (2003), no. 3, 401--408.
 
 \bibitem{ka3} N.J. Kalton, {\it The three space problem for locally bounded F-spaces}, Compo. Math. 37 (1978), 243--276.
 
 \bibitem{KM} N.J. Kalton and S. Montgomery-Smith, {\it Interpolation of Banach spaces}, Handbook of Geometry of Banach Spaces, Vol. 2, (W.B. Johnson and J. Lindenstrauss, editors), Elsevier, Amsterdam, 2003, 1131--1175.

\bibitem{KaPe} N.J. Kalton and N.T. Peck, {\it Twisted sums of sequence spaces and the three-space problem}, Trans. Amer. Math. Soc. 255 (1979), 1--30.



%
%



%
%
%
%
%
%
%
%
%

%
\bibitem{LinTz} J. Lindenstrauss and L. Tzafriri, \emph{Classical Banach
spaces I, sequence spaces}. Ergebnisse der Math und ihrer Grenzgebiete, Vol. 92. Springer-Verlag, Berlin-New York, 1977.
%


\bibitem{Mau2} B. Maurey, {\it A remark about distortion}, Geometric aspects of functional analysis (Israel, 1992–1994), 131–142, Oper. Theory Adv. Appl., 77, Birkha\"user, Basel, 1995.

\bibitem{Mau3} B. Maurey, {\it Type, cotype and K-convexity}, Handbook of the geometry of Banach spaces, Vol. 2, 1299–1332, North-Holland, Amsterdam, 2003.

%


\bibitem{JS3} D. Morales and J. Su\'arez de la Fuente, {\it Some more twisted Hilbert spaces}, Annales Fennici Mathematici, 46(2), 819--837.

\bibitem{NT} N.J. Nielsen and N. Tomczak-Jaegermann, {\it Banach lattices with property (H) and weak Hilbert spaces}, Illinois J. Math. Volume 36, Issue 3 (1992), 345--371.

\bibitem{OS} E. Odell and T. Schlumprecht, {\it The distortion problem}, Acta Math. 173 (1994), no. 2, 259–-281.

\bibitem{OS2} E. Odell and T. Schlumprecht, {\it Distortion and asymptotic structure}, Handbook of the geometry of Banach spaces, Vol. 2, 1333-–1360, North-Holland, Amsterdam, 2003.
\bibitem{Pi} G. Pisier, {\it Weak Hilbert spaces}, Proc. London Math. Soc. (3) 56 (1988), no. 3, 547--579.

\bibitem{Schlum} T. Schlumprecht, {\it An arbitrarily distortable Banach space}, Israel J. Math. 76 (1991), no. 1-2, 81–-95.

%

\bibitem{JS} J. Su\'arez de la Fuente, {\it A weak Hilbert space that is a twisted Hilbert space}, J. Inst. Math. Jussieu 19 (2020), no. 3, 855--867.

\bibitem{JS2} J. Su\'arez de la Fuente, {\it A space with no unconditional basis that satisfies the Johnson-Lindenstrauss lemma},  Results Math. 74 (2019), no. 3, Art. 126, 14 pp.

\bibitem{JS4} J. Su\'arez de la Fuente, {\it A universal formula for derivation maps and applications}, Analysis Mathematica (to appear).

\bibitem{ToJ} N. Tomczak-Jaegermann, {\it Distortions on Schatten classes $C_p$}, Geometric aspects of functional analysis (Israel, 1992–1994), 327–334, Oper. Theory Adv. Appl., 77, Birkh\"auser, Basel, 1995. 

\bibitem{Wat} F. Watbled, {\it Complex interpolation of a Banach space with its dual}, Math. Scand. 87 (2000), no. 2, 200--210. 




\end{thebibliography}
\end{document}